# Locating Platforms and Scheduling a Fleet of Drones for Emergency Delivery of Perishable Items

Monica Gentili [a], Pitu B. Mirchandani [b], Alessandro Agnetis [c], Zabih Ghelichi [a]

[a] Department of Industrial Engineering, University of Louisville, KY, USA
[b] School of Computing, Informatics, Decisions System Engineering, Arizona State University, AZ, USA
[c] Department of Information Engineering and Mathematical Science, University of Siena, Siena, Italy

**Abstract**
Motivated by issues dealing with delivery of emergency medical products during humanitarian disasters, this paper addresses the general problem of delivering perishable items to remote demands accessible only by helicopters or drones. Each drone operates out of platforms that may be moved when not in use and each drone has a limited delivery range to service a demand point. Associated with each demand point is a disutility function, or a cost function, with respect to time that reflects preferred delivery clock time for the demanded item, as well as the item's perishability characteristic that models nonincreasing quality with time. The paper first addresses the problem of locating the platforms as well concurrently determining which platform serves which demand points and in what order – to minimize total disutility for product delivery. The second scenario addresses the two-period problem where the platforms can be relocated, using useable road network, after the first period. It can be easily proven that continuous time versions of these problems are NP-Hard. However, a practical "time-slot" version of the problem, where time is discretized into slots, can be solved by standard optimization software. Extensive computational experiments, using different drone delivery ranges as well as different drone fleet sizes, provide valuable insights on the performance of such drone delivery systems.

**Keywords:**
Drones, Delivery of Perishable Items, Time-slot scheduling and routing. Emergency Delivery Systems, Humanitarian Logistics.

## 1 Introduction

Road-network based emergency public safety systems, such as ambulance, fire protection, and police with associated traffic management protocols and logistics support systems have been designed and developed over the last half century. There are also complementary airborne emergency systems, primarily using helicopters, that system engineers and operations researchers have been developing for emergency humanitarian assistance, for example search-and-rescue missions [1, 2], delivery of emergency goods [3-5], and emergency medical assistance [6]. Given the rapid developments in technologies for remote controlled unmanned aerial vehicles (UAVs - also referred to as unmanned aircraft systems, UAS, and *drones*), this paper argues that drones and associated logistical systems can cost-effectively address some problems of emergency delivery of medicines to patients that are remote and not easily accessible via roads from medical facilities and inventories.

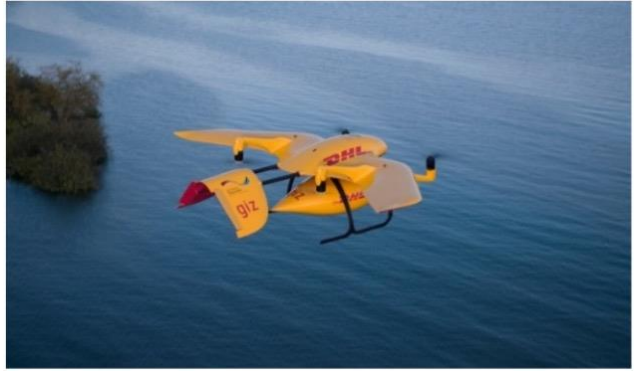

**Figure 1. Illustration of the Drone Delivery Technology: DHL Parcel copter.**
https://www.dpdhl.com/en/media-relations/specials/dhl-parcelcopter.html

Advantages of drones for the delivery of goods from anywhere to anywhere (illustrated in Figure 1) include (1) vertical takeoff from and landing in small places, (2) low capital and operating costs, (3) no requirements of a pilot on board, and (4) no need for a costly expansive infrastructure for takeoff/landing.

Uses of such drones are rapidly increasing. Recently, operational engineers and researchers have been developing logistical tools for operating drones for a wide range of applications. Some recently reported applications and attendant logistical models are for traffic surveillance [7, 8]; parcel delivery [9]; traffic planning and management [10]; agriculture surveillance and crop spraying, monitoring fires, sports and entertainment event coverage [11] among many others.

This paper considers the problem of delivering emergency medical supplies, although the developed logistical models can be used for other applications where low payloads are to be transported or delivered or picked up quickly to/from points that are not easily and quickly reachable by roads. In particular, this paper addresses the problem of sending perishable medical supplies to demand points that are not well connected to the rest of the network, for example, communities which are accessible only via damaged or poor quality roads during natural disasters such as hurricanes and floods.

Describing our problem more specifically, we will assume there are $p$ delivery drones to serve $m$ remote demand points that need delivery of emergency medical supplies such as blood units. These demand points are reachable only by drones where each of the drones is constrained by a limited distance range to service a demand point. The drones operate out of mobile platforms (see Figure 2 for an example) which may be moved on traversable roads. Each demand point requires a single package of product whose utility decreases with delivery time due to product perishability. Furthermore, a customer's demand may have a delivery time requirement. The main problem addressed is to locate $p$ platforms, and their associated drones, so that the total disutility for delivery times is minimized.

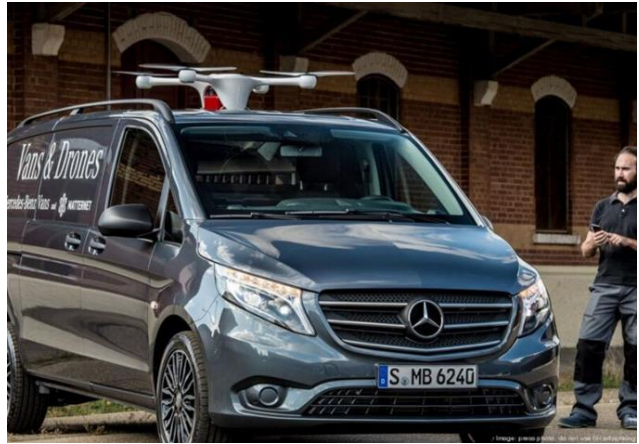
**Figure 2. Matternet, a logistics company in Switzerland, received official approval in March 2017.**
*https://www.techspot.com/news/71185-self-flying-drone-network-completes-maiden-delivery-switzerland.html*

This paper is motivated by the recent 2017 hurricanes in USA, Harvey and Irma. It has been observed by both government agencies and organizations that the delivery of medical supplies for patients in inundated regions is a unique and challenging problem. For example, Forbes.com [12] noted in August 2017, right after Hurricane Harvey, that subsequent to power restoration, "delivering medical care to individuals in compromised living quarters – whether in a flooded home or in a shelter – presents unique challenges to first responders, as well as doctors and nurses who serve on the front lines". Right after Hurricane Irma, in September 2017, Aircraft Owners and Pilots Association (AOPA) volunteered to deliver medicine with piloted aircrafts to flooded areas in Florida and Virgin Islands due to Hurricane Irma [13]. Such observations motivate the development of a drone-based emergency medicine delivery systems, which this paper proposes.

The paper addresses two variants of the problem. In the first variant, we assume there is a fixed time, for example one day, from morning to evening, to deliver the products demanded. In the second variant, we consider two periods, say two days, to deliver the demanded products. In both these problem variants, the objective is to minimize the total disutility for the delivery of the required packages of products. In the second variant, it is assumed that the mobile platforms may be moved in the evening of the first day for second day operations. There are other extensions of these problems that we discuss in the concluding section.

Previewing the paper, the next section discusses some related reported research. Section 3 formulates the optimization problem for the first variant of the problem, namely the single-period variant, and reports extensive computational experience for a case study. Section 4 formulates the two-period variant of the problem and discusses some computational experiences. The concluding section 5 discusses some extensions of the problem and possible directions for future research.

Briefly summarizing the contributions of the paper, which will become clearer shortly, are as follows:
- It addresses a new application of UAVs for the delivery of aid packages to disaster-affected areas. It extensively studies different operational aspects of a drone delivery system and derives multiple practical managerial insights, through computational analyses of various cases.

- It introduces a new "timeslot" formulation to address the integrated drone location and scheduling problem which gives us the opportunity to efficiently solve the problem by using off-the-shelf solvers in reasonable amounts of time.
- It adopts a very general objective function that, combined with the time slot formulation, captures the cost or disutility at delivery times for each individual demand, as long at cost/disutility is non-decreasing in time.

## 2 Related literature

### 2.1 Some implemented and proposed applications

A recent survey by Otto et al. [14] provides a comprehensive list of actual and potential civil applications of drones, from search and find, data gathering, measuring network performance, surveillance, product delivery, to emergency and medical applications. Search and find operations have a long history in field operation of drones; they have been used for wildlife monitoring, for search and rescue applications after disasters, for search and destroy for military operations, and others. Drones have been very useful for data gathering on the status of physical infrastructures (e.g., tracking the progress of construction sites, assessment of inventories, maintenance and inspections of bridges and roads, evaluation of damages and risks after natural disasters, etc.), for searching and finding lost people or targets, and for agricultural applications for monitoring the status of growth, irrigation, etc.

Drones have also been proposed for establishing wireless communication. Multiple drone base stations can provide wireless communication network to serve users [15]. In the event of a natural disaster or simply for improving the communication/internet performance, drones and balloons are being planned as ad-hoc communication nodes by Google's project Loon [16] and Facebook's Aquila project [17].

Recently, drones have enabled organizations to transport packages of different sizes and weights, for example, they have been used for postal service to a rural region of Southern France [18]. In a recent paper, drones were proposed for resupplying delivery trucks for same-day delivery systems [19]. In the context of delivery of health-related items, they have been used for blood transport in Rwanda [20], tested for delivery of health equipment such as defibrillators [21] and surgical robots [22] to remote areas. Rosser et al. [23] provide a review on the use of drones in healthcare.

### 2.2 Logistics issues and some reported results

The logistics issues related to a fleet of drones consider the types of application, the time or frequency of each service or delivery, the number of drones in the fleet (we will denote this number by notation $p$), the transport characteristics (e.g., direct straight-line delivery vs. tour of deliveries), and the performance objective for the fleet of drones. Among the first uses of drones was to search an area $R$ and find a "target" in $R$. Typically, in such problems we often have a gridded map of $R$ or some sort of partitioned map where each cell has its own probability of finding the target; in case of a moving target, the probability of the cells evolves over time [24-26]. The related logistics problems attempt to find paths or cover regions to maximize the probability of finding the target (see e.g., [27]).In monitoring traffic incidents or surveying disaster regions, drones can be sent to visit locations within specified time windows [28-31] to optimize a path or tour-cost objective.

In the context of facility location problems for drone delivery systems, Chauhan et al. [32] propose a capacitated facility location problem where drones perform one-to-one trips between facilities and demand locations. The problem of determining optimal locations for recharging stations for commercial drone delivery systems was studied by Hong et al. [33]. Dukkanci et al. [34] developed a nonlinear model for a drone delivery problem to find optimal launching locations and the speeds at which drones travel. They assume drones are transported to launching locations by traditional vehicles. In a recent study, Salama and Srinivas [35] proposed a truck-drone delivery model for last-mile delivery, where the drones deliver from a single movable truck. The customers are clustered around focal points while the truck travels from focal point to focal point. They optimized several decisions including clustering of customers, the location of focal points, and truck routes to deliver each customer by either a drone or the truck aimed at minimizing total delivery costs and/or delivery makespan.

Given the challenges of disaster relief management and the importance of timely emergency response, the application of drones for humanitarian logistics has recently received some attention among researchers. Kim et al. [36] studied a platform location problem where the operating conditions due to weather, fuel consumption and flight distances introduce uncertainties that were handled by chance constraints; the problem was modeled as minimizing the needed drones to cover the demands where each demand requires a single drone delivery. Rabat et al. [25] studied the problem of delivering light-weight items from a single depot to demand points with three levels of priorities and drones could be refueled at specified locations. The proposed model assumes single delivery at each demand point and considers energy consumption and payload weight to minimize the total cost of routing of the drone fleet for serving the demands. Recently, Zubin et al. [37] used a standard Vehicle Routing Problem formulation, with a solution approach that uses a large-scale neighborhood search heuristic, for home delivery of pharmaceutical items using drones; they conducted an extensive evaluation based on a realistic case study.

As noted by Otto et al. [14], a large number of papers have dealt with routing problems to a set of locations, giving rise to models which are typically related to the Traveling Salesman Problem (TSP) or the Vehicle Routing Problem (VRP) [38]. In many studies, payload limitations imply that a drone can carry only one package at a time [39, 40]. In such cases, the situation is not adequately captured by a routing model, since the drone performs a separate flight for each destination, always returning to the base between two flights. In this case, the key decision concerns the sequence in which different destinations are served.

In an increasing number of studies, drone operations are coordinated with those of other road vehicles or marine vessels. The roles of the different vehicles depend on the specific scenario being considered. Usually, either the performance of the drones or of the other vehicles may be the main logistic objective [14]. This paper is concerned with the scenario where vehicles (such as trucks) support drones operations, that is where trucks constitute the platforms from which a drone starts and ends each flight, and where drones are refueled, or their spent batteries are swapped with charged ones [40, 41]. Unlike cases where drones and vehicles independently perform delivery tasks [4, 9], in this scenario, the trucks only play the role of movable platforms, and do not perform any delivery because the roadways have limited accessibility and destinations are remotely located where drones become essential to ensure timely delivery. Recently, results have appeared where drones work in cooperation with other vehicles (trucks, cars, trains, and sea vessels) to meet surveillance or delivery objectives [42-44]. The underlying models for such problems are often variants of TSPs and VRPs. This paper does not address the TSP/VRP variants that may result when the drone platforms on trucks and drones themselves are both delivering emergency medicine packages; this is left for potential future research as discussed in Section 5.

**2.2 Contributions of this research**

Following are the contributions of this research, both methodological and practical:

1. First, we consider a new application of UAVs to address the problem of timely delivery of time-sensitive aid items to disaster-affected areas. We develop multiple optimization models to study different operational aspects of a drone delivery system through an extensive set of experiments. The operational aspects that we studied in this research include (i) the effect of using homogenous fleet and heterogeneous fleet of drones (short-range vs long-range drones), (ii) the effect of variants of disutility/cost functions (further discussed in point 3 below) , (iii) the effect of discretization of time, (iv) the effect of number of drones to utilize per platform, and (v) allowing an additional period during which platforms can be relocated. Based on these experiments, we gained some interesting operational insights on the performance of the system to facilitate decision-making in humanitarian logistics.

2. Second, we develop optimization models for **concurrently** locating drone platforms and scheduling a set of trips for each drone. We first develop a classical Integer Linear Programming (ILP) model (see Appendix) to formulate a single period drone location and scheduling problem. However, we prove that the drone location and scheduling problem is *NP-hard* in strong sense. We then introduce a novel space-time formulation which discretizes the time horizon into multiple timeslots that turns out to have several practical advantages. First, this new approach gives us the opportunity to efficiently solve the drone location and scheduling problem using off-the-shelf solvers, e.g., CPLEX, in reasonable amounts of time. Second, the granularity of the timeslots determines the level of approximation such that the smaller the timeslots, the more precise the solutions. That is, there is a trade-off between the level of approximation and computational complexity of the model. Thus, the developed formulation is highly flexible in terms of determining the level of approximation and controlling the computational effort (see Section 3.2.2.2). Third, this formulation can accommodate the case where a heterogeneous fleet of drones (short-range and long-range drones) can be utilized (see Appendix). The formulation can also be extended to consider a two-period model where the drone platforms can be relocated after the first period.

Last, but not least, a distinctive feature of our model is its ability to model a variety of optimization objective functions, not necessarily linear or continuous, due to the flexibility of the time-slot formulation to capture individual demand's cost. In most models where drones are to be used for goods delivery, the predominant performance indicators are typically system-oriented, reflecting efficiency or productivity requirements for the system. Such objectives include minimizing the time for completing all delivery operations (e.g., [45]) or minimizing operational costs, typically corresponding to the cost of a tour in a routing application [46, 47]. Rather, our model focuses on individual recipients' objectives and needs. Indeed, in the delivery of perishable items, depending on specific features of the product being delivered, the utility of a recipient may vary considerably over time, and such utility may be very different for recipients at various demand points. For instance, the quality of a blood unit can be acceptable up to a given time, after which the product starts deteriorating at an increasing rate, eventually becoming unacceptable. Furthermore, the urgency of each delivery is dictated by the need of providing timely relief (one may refer to it as "due time" for delivery) for the corresponding patient and therefore a recipient's utility may be different given which other patients are to be serviced by the assigned drone. The common feature of these disutility/perishability functions is that they can be highly non-linear and non-decreasing with respect to time of delivery. Linear increasing cost function is the special case that is used in most delivery models. In the model developed, our time-slot-formulation can easily handle any non-decreasing function, where for each demand point $k$, a *disutility* $f_k(t)$ is specified expressing the disutility of delivering a package at point $k$ at time-slot representing time $t$. We do not make any assumption on the shape of $f_k(t)$, except that it is non-decreasing as $t$ grows since a late delivery is not better than an earlier delivery.

## 3 The Single-period Problem

Consider first the single-period problem. There is a set $P$ of $m$ candidate platform sites and a set $D$ of $n$ demand points. A total of $p \leq m$ platforms must be located. Each platform hosts one drone. If a platform is at site $i$, a drone departing from it can serve the demand points located within a specified distance from it, denoted by $D(i) \subseteq D$. Each demand point $k$ has a request for a package that is preferred to be delivered within a given due date $d_k$, and which can be shipped by any platform at site $i$ such that $k \in D(i)$. A drone can only carry out one package at a time, and it must head back to the platform after carrying out each delivery. We let $p_{ik}$ denote the roundtrip time from platform at site $i$ to point $k$ (inclusive of the time needed to unload the package and send the drone back). We will assume that if a drone leaves platform at site $i$ at time $t$ for servicing demand point $k$, the demand point is reached at time $t + \frac{p_{ik}}{2}$. A given setup time (sequence-independent) must occur between two consecutive flights of the same drone, to allow for refueling or battery swap and securing a new package to the drone.

For each demand point $k$, a *disutility function* $f_k(t)$ is specified expressing the disutility or cost of delivering the package of point $k$ at time $t$. We do not make any assumption on the shape of $f_k(t)$, except that it is non-decreasing as $t$ grows.

The problem is to concurrently decide (a) where to locate the platforms, (b) which demand points must be served from each platform, and (c) in which sequence these demand points are to be served, so that all packages are delivered within a given time period $T$ (e.g., one day) and the total disutility is minimized.

Summarizing, we consider the following *Drone Location and Scheduling problem* (*DLS*):

*Given $m$ candidate platform sites and $n$ demand points, each demand point $k$ having its individual disutility function $f_k(t)$, given the roundtrip time $p_{ik}$ from site $i$ to demand point $k$ for all $i$, and for all $k$, find the location of $p$ launch platforms and the optimal schedule of drone flights so that all the demand points are served within a period $T$ and the total disutility is minimized.*

We note that, once the platform locations are selected, DLS reduces to a pure *unrelated machine scheduling problem* [48], in which jobs correspond to the demand points and processing times relate to distance between pairs of demand points and platforms. To the best of our knowledge, DLS has never been addressed before.

Let us now consider the complexity of DLS. The following result can be easily established.

**Theorem 1**. *DLS is NP-hard in the strong sense.*

*Proof.* Consider an instance of the parallel-machine scheduling problem $Pm||\sum w_k C_k$, in which $n$ jobs of integer processing times $p_k$ have to be allocated to $m$ machines to minimize total weighted completion time. The problem consists of deciding whether a schedule exists having total weighted completion time (*TWCT*) not exceeding an integer $K$. We can define an instance of DLS in which: (*i*) $p = m$, (*ii*) all $m$ sites coincide (so that the roundtrip time to reach any demand point is the same for all sites), (*iii*) for each job in the instance of $Pm||\sum w_k C_k$ there is a demand point $k$ such that the processing time (i.e., the roundtrip time) is precisely $p_k$, (*iv*) the time horizon $T$ is long enough not to be binding (e.g. $T = \sum p_k$), and (*v*) disutility is defined as $f_k(t) = w_k t$. We want to establish whether a solution exists having total disutility not exceeding $K - (\frac{\sum_k w_k p_k}{2})$. □

Consider a schedule for $Pm||\sum w_k C_k$, where $C_k$ denotes the completion time of job $k$, so that its total weighted completion time is $TWCT = \sum_k w_k C_k$. We can build a solution to DLS identifying the $m$ platforms with the $m$ machines, along with the corresponding schedules. Since demand point $k$ is reached exactly after $\frac{p_k}{2}$ time units from its departure, the total disutility of such a solution is given by $\sum_k w_k \left(C_k - \frac{p_k}{2}\right) = TWCT - (\sum_k \frac{w_k p_k}{2})$. Hence, the instance of DLS is a yes-instance if and only if the instance of $Pm||\sum w_k C_k$ is a yes-instance.

### 3.1 Time-slot formulation for the single-period DLS

Our first approach was to model the single-period DLS as a classical ILP model (which is provided in the Appendix), using assignment and sequencing variables. Some preliminary computational results showed that such a formulation is not particularly efficient, and only small instances could be solved to optimality within an acceptable computational time. We therefore developed a more useful and efficient ILP formulation, described in this section, which exploits a key issue of DLS: flight durations may be expressed in terms of discrete time units.

We examined the literature to come up with reasonable characterization of drones with respect to range and speeds. Currently, a large variety of drones exist and new ones are being developed for non-military applications (see [49, 50]). A maximum range of 30-50 km and a nominal speed of 60 km/h are assumed in our computational experiments for the emergency delivery system being modeled and analyzed here[1]. Also, we envision that flight times, and delivery and refueling/recharging durations are normally less than 15 minutes. Hence, in our scenarios we assumed 15 minutes duration as our discrete time unit, and values for the round-trip time $p_{ik}$ will take values between one and eight time units (i.e., 15 min to 2 hours). Although in principle a delivery time may require a non-integer number of such time units, we assuming that rounding up these times to integers should be acceptable; in the worst case it might result in a conservative schedule. We emphasize that these ranges, speeds and time units are assumed only for our experiments; they can easily be replaced by actual values in any field implementation.

Under these assumptions, we discretize the time horizon in a predefined number of time slots whose index set is denoted by $T$. So adopting, for example, 15-minute time slots, a one-day period of 12 hours is equivalent to 48 time slots, i.e., $T = \{1,2,\ldots,48\}$. The resulting formulation requires only two variable types:
- $y_i \in \{0,1\}$ which is equal to one if site $i$ is selected, and equal to zero otherwise.
- $x_{itk} \in \{0,1\}$ which is equal to one if a drone gets back to platform at site $i$ at the end of time slot $t$ after having served demand point $k$.

In this formulation, time spans (such as roundtrip times $p_{ik}$) are always expressed by a number of time slots. So, $f_k(t)$ denotes the disutility incurred by $k$ if it is served at the end of time slot $t$. A non-integer $t$ indicates that a drone may be served in the middle of a time slot. So, if a drone *is back at platform at site $i$ at time slot $t$ after* serving demand point $k$ (i.e., if $x_{itk} = 1$), the disutility incurred by $k$ is $f_k(t - \frac{p_{ik}}{2})$. Note that no assumption is made on the specific shape of the disutility function. The problem can be formulated as follows.

$$\min \sum_{i \in P} \sum_{k \in D(i)} \sum_{t \in T} f_k\left(t - \frac{p_{ik}}{2}\right) x_{itk}$$

$$\sum_{i \in P} \sum_{t \in T} x_{itk} = 1 \qquad \forall k \in D \qquad (1)$$

$$\sum_{\substack{h \in D(i) \\ h \neq k}} \sum_{\tau = t - p_{ik}}^{t} x_{i\tau h} \leq M_{ik}(1 - x_{itk}) \qquad \forall i \in P, \forall t \in T, \forall k \in D(i) \qquad (2)$$

$$x_{itk} \leq y_i \qquad \forall i \in P, \forall t \in T, \forall k \in D \qquad (3)$$

$$\sum_{i \in P} y_i = p \qquad (4)$$

$$\sum_{k \in D(i)} \sum_{t \in T: t \leq p_{ik}} x_{itk} = 0 \qquad \forall i \in P \qquad (5)$$

$$x_{itk} \in \{0,1\} \qquad \forall i \in P, \forall t \in T, \forall k \in D \qquad (6)$$

$$y_i \in \{0,1\} \qquad \forall i \in P \qquad (7)$$

---

[1] For example, we may envision the use of a DHL Parcelcopter 4.0, which is a tiltrotor with maximum payload of 4 kg, maximum range of 65 km and maximum speed of 130km/h: https://www.dpdhl.com/en/media-relations/specials/dhl-parcelcopter.html

Constraints (1) impose that each demand point is served exactly once. Constraints (2) are capacity constraints, implying that whenever a drone is flying towards or back from $k$, it cannot carry out other missions. Parameter $M_{ik}$ is a number expressing an upper bound on the number of different trips a drone can make during a time span of length $p_{ik}$. If $x_{itk} = 1$, these constraints imply that no other point $h \neq k$ can be served in the time interval $[t - p_{ik}, t]$ from platform at site $i$. Constraints (3) ensure a demand point is served from a drone departing from platform at site $i$ only if platform at site $i$ is selected. Constraint (4) requires the selection of exactly $p$ platform sites. Constraints (5) account for the duration of the first flight of the day from each platform site $i$.

Note that, this formulation allows one time slot between two consecutive delivery missions, to account for setups (refueling/recharging, loading/unloading etc.). The model can be easily adjusted if such setup times are supposed negligible (as in [26]) or require more than one time slot. Additionally, the formulation assumes a homogenous fleet of drones (i.e., drones of the same type) is used. The mathematical formulation can be modified to account for a mixed fleet of drones. In particular, in the Appendix (see formulation (27) – (38)) we provide a modification of model (1)-(7) where a mixed fleet that includes two types of drones (short-range and long-range drones) is considered. Specifically, the modified formulation assumes a fleet of $N = N^{(1)} + N^{(2)}$ drones is available, where $N^{(1)}$ is the number of short-range available drones, and $N^{(2)}$ is the number of long-range available drones.

Finally, note that formulation (1) - (7) (as well as formulation (27)-(38) in the Appendix) can be used for the case where more than one drone can be allowed on the same platform. In order to do that, no modification of the structure of the mathematical model is needed. This will be better explained in subsection 3.2.2.3 where we consider the case where up to two drones can be used on each platform.

**3.2 Computational experiments for the single-period problem**

In this section, we report on computational experiments designed to test the feasibility of solving real-size problems with the time-slot formulation (1) - (7). We solved the model using CPLEX 12.8 on a server with CPU Intel Xeon E5-2650 v3 @ 2.30GHz, with 128 GB of RAM. A maximum time limit of three hours was set for each run.

*3.2.1 Description of the computational experiments*

Our experiments consider a reference scenario, namely Central Florida. Over the road network, we identified a set of candidate platform sites, while the demand points were randomly scattered away from the main roads. Data common to all experiments are:

- Number of demand points: 100
- Number of candidate platform sites: 25
- Time horizon for the single-period model: 12 hours (8am to 8pm)

Figure 3 shows the reference scenario. Squares denote potential platform sites, which are located along main roads and freeways, while circles indicate demand points, which are scattered throughout rural areas.

While the above parameters are common to all experiments, we have considered different scenario settings to evaluate how the solution of the model changes with different parameter values. Specifically, we considered:

- Maximum flying range of a drone ($T_{max}$): 30km (short-range drones) and 50km (long-range drones).
- Number of platforms: $p = 1, 2, \ldots, 25$.
- Time slot length: 5 minutes or 15 minutes.
- Type of fleet: either homogenous (i.e., only short-range drones or only long-range drones) or mixed fleet type.
- Number of drones allowed per platform: either one drone per platform or up to two drones per platform.
- Shape of disutility function (see below).

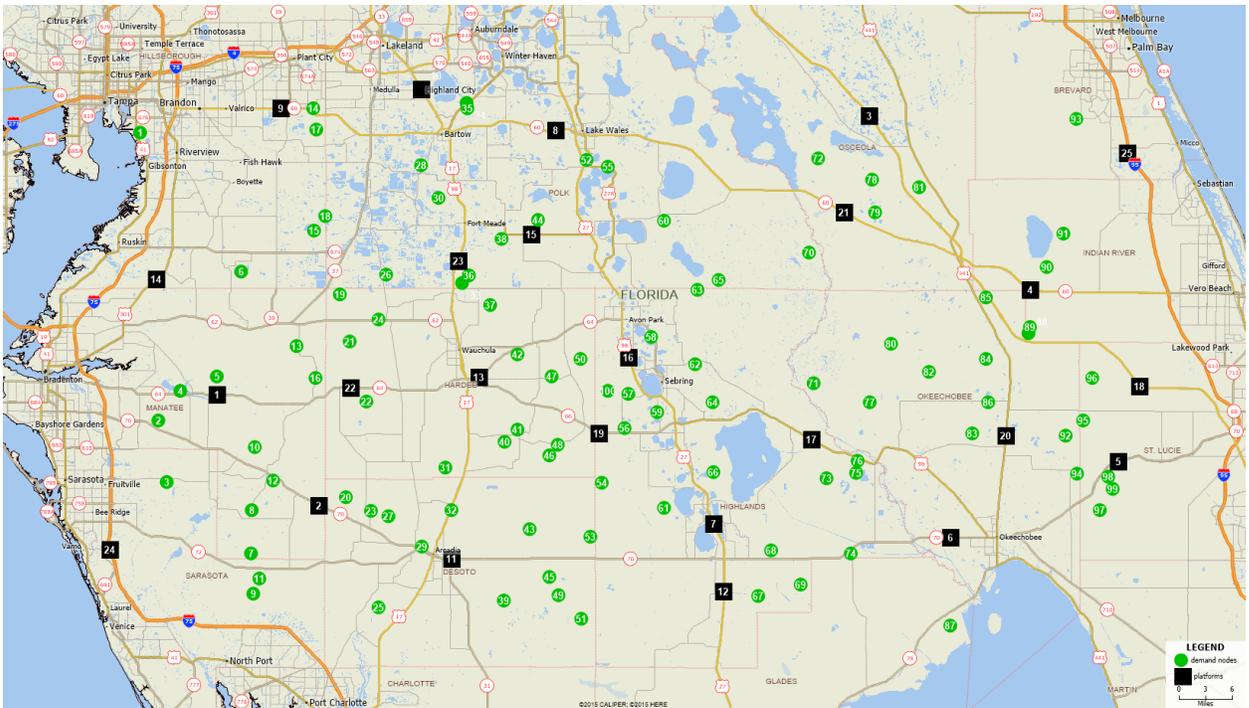

**Figure 3. Map of central Florida, with 25 potential platform sites (squares) and set of 100 demand points (circles).**

As discussed in Section 1, a specificity of our setting is that two distinct factors could contribute to the disutility of the delivery time. The first factor expresses perishability, that is, the quality of the item decreases over time, as typical of certain health-related products such as plasma or other blood derivatives. The second factor expresses the fact that, from the viewpoint of the customer, disutility increases after a certain due date $d_k$, because, for example, if no treatment is received within $d_k$, the clinical conditions of the recipient may start deteriorating. Although it is not necessary that these factors be given as functions of delivery time $t$ (that is, they could be table look-up from empirical data on perishability and demand priorities), for the computational experiments on a case study we assumed that the disutility functions $f_k(t)$ to be the sum of two components, one for perishability and one for due date, as described below by equations (8)-(10), both of which increase with $t$.

$$f_k^{(1)}(t) = \frac{A_k}{100^2} t^2 \tag{8}$$

$$f_k^{(2)}(t) = \begin{cases} 0, & t \in [0, d_k] \\ \frac{B_k}{(100 - d_k)^2} (t - d_k)^2, & otherwise \end{cases} \tag{9}$$

$$f_k(t) = f_k^{(1)}(t) + f_k^{(2)}(t) \tag{10}$$

where $A_k$ is a perishability coefficient for service to demand point $k$, while $B_k$ and $d_k$ are importance coefficients and the due date associated with that demand, respectively. Adjusting these coefficients easily allows us to study effects of perishability and importance in these logistics problems. We note that there is a large body of literature in social sciences and operations research ondeveloping utility functions (or cost functions or value functions) based on stated or revealed preferences, see for examples [51-56]; some empirical efforts must be expended to develop such functions or simple tables, for actual application of the models developed here.

In our experiments we considered two case scenarios. For the general case, the *non-uniform scenario*, demand points may have different disutility functions, that is, individual due dates and perishability coefficients are specified for each demand point, to reflect the specificity of requirements for each shipment. In the computational experiments, $A_k$ and $B_k$ for each demand point $k$ are randomly generated between 50 and 1000, and $d_k$ is randomly chosen between 2 and 8 hours. For comparison purposes, a second case, the *uniform scenario,* is studied where timeliness and perishability do not depend on the specific demand point, that is, where all demand points have the same due date $d$ and the same perishability and importance coefficients. In the computational experiments for this uniform scenario we set $A = B = 100$ and $d = 2$ hours.

Before examining the computational results, we illustrate a sample solution. Figure 4 shows the solution of the single-period model obtained for $T_{max} = 50$, $p = 9$, and time slot interval equal to 15 minutes. Black squares indicate unselected platform sites; red stars indicate selected platform sites. Each selected site is linked to the allocated demand points (circles in green). The number in each green circle is the time slot in which the drone is back to the platform after visiting the demand point, so it indicates flight schedule information for each drone. To illustrate, the drone from selected site #1 first visits the demand point just west of it and returns at time $t = 2$, then delivers at the demand point just north of the platform and return at time $t = 4$, then visits the demand point southeast of the platform and returns at time $t = 7$, and so on, the order of the deliveries being inferred by the numbers in the green circles.

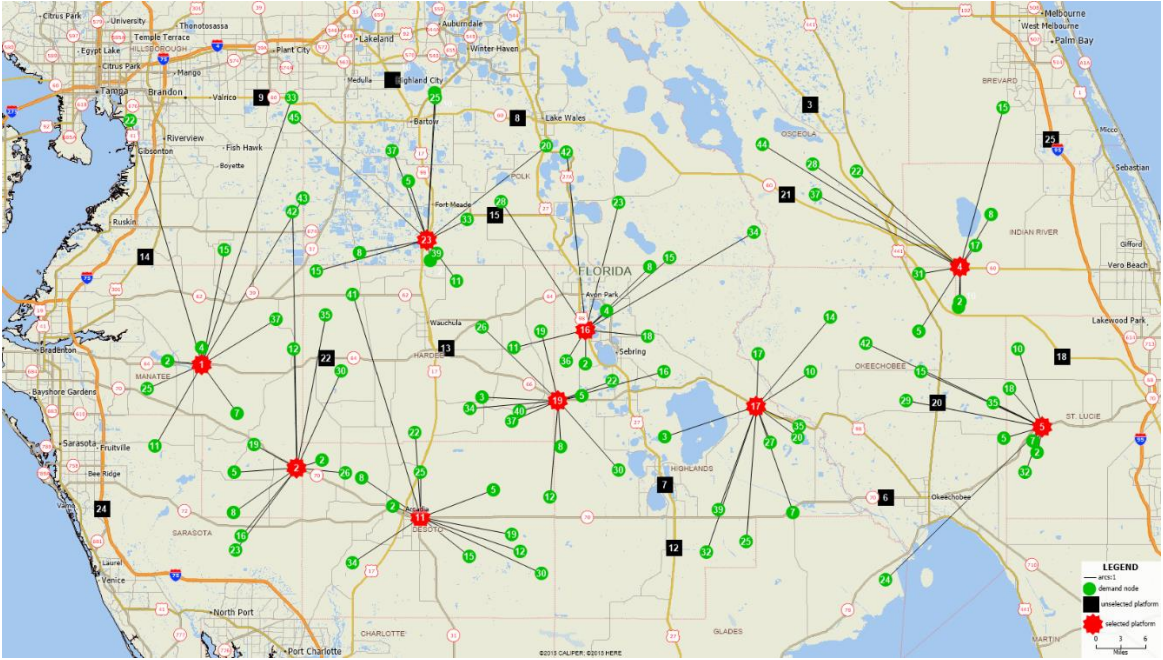

**Figure 4. Example solution for the single-period model.**

*3.2.2 Computational results*

The results are organized into the following subsections:
- In subsection 3.2.2.1, we analyze the effects of using different homogeneous fleet of drones (i.e., short-range vs long-range drones) and of having uniform disutility functions vs non-uniform disutility functions.
- In subsection 3.2.2.2, we analyze the effect of using different values of time slot length.
- In subsection 3.2.2.3, we analyze the effect of allowing up to two drones per platform.
- In subsection 3.2.2.4, we analyze the effect of using a homogenous fleet of drones versus using a mixed fleet of drones.

The input data and the detailed results are available upon request to the authors.

*3.2.2.1 Short-range vs long-range drones, and uniform vs non-uniform disutility functions.*

Results in this section are obtained by solving model (1) - (7) with 15-minute time slots. Results are summarized in Tables 1-5. Several insights can be gained from examining the results, as discussed next.

*Observation 1: For both uniform and non-uniform scenarios, if the number of drones/platforms is too small no allocation and schedule can service all demand points within one day.*

Table 1 shows the computational time and the disutility function value for each value of $p$ and $T_{max}$ in the non-uniform scenario; Table 2 give these results for the uniform scenario. The tables also show the marginal improvement (in %) of the disutility value with respect to the values corresponding to the minimum feasible $p$ for that scenario. For both scenarios no feasible solution can be found for $p \leq 10$ when $T_{max} = 30$, and for $p \leq 8$ when $T_{max} = 50$. All the instances were solved within the three hours limit set for computation time, with the exception of $T_{max} = 50$, $p = 8$ for the non-uniform scenario, and $T_{max} = 50$, $p = 7, 8, 9, 10$ for the uniform scenario for which it was not possible to certify infeasibility nor to find a feasible solution within the time limit.

*Observation 2: For both uniform and non-uniform scenarios, the marginal benefit of having one more platform monotonically strictly decreases.*

For the non-uniform scenario (Table 1), when going from $p = 11$ to $p = 25$ with $T_{max} = 30$, the disutility drops from approximately 1,185 to 185, but half of the benefit is due to the first three additional drones. Similarly, when $T_{max} = 50$, going from $p = 9$ to $p = 25$ the disutility monotonically decreases, from 1,699 to 183. Increasing $p$ from 9 to 14 yields almost 70% of the saving, while the addition of 11 more platforms only contributes the remaining 30% improvement. Similar observations can be made for the uniform scenario (Table 2). The disutility values cannot be directly compared with the non-uniform scenario, since the instances were generated differently. However, observe that the relative disutility decrease attained when $p = 25$ with respect to the minimum feasible $p$, is slightly larger in the uniform case (88.7% for $T_{max} = 30$ and 92.2% for $T_{max} = 50$, respectively) than in the non-uniform case (84.4% for $T_{max} = 30$ and 89.3% for $T_{max} = 50$). Also, computational times are larger in the uniform case, perhaps because of the presence of nearly equivalent solutions.

Table 1. Computational time (sec.) and disutility function values for the non-uniform scenario.

| | Non-Uniform Scenario | | | | | |
|---|---|---|---|---|---|---|
| | $T_{max} = 30$ | | | $T_{max} = 50$ | | |
| $p$ | CPU time (sec.) | Disutility | Marginal Benefit | CPU time (sec.) | Disutility | Marginal Benefit |
| 1 | 50.43 | Infeasible | - | 141.11 | Infeasible | - |
| 2 | 49.22 | Infeasible | - | 196.89 | Infeasible | - |
| 3 | 56.13 | Infeasible | - | 291.28 | Infeasible | - |
| 4 | 62.00 | Infeasible | - | 345.60 | Infeasible | - |
| 5 | 63.50 | Infeasible | - | 1,273.60 | Infeasible | - |
| 6 | 77.59 | Infeasible | - | 2,100.41 | Infeasible | - |
| 7 | 75.73 | Infeasible | - | 2,574.14 | Infeasible | - |
| 8 | 96.18 | Infeasible | - | 10,808.22 | Infeasible | - |
| 9 | 88.88 | Infeasible | - | 10,811.92 | 1,699.67 | - |
| 10 | 129.09 | Infeasible | - | 6,501.91 | 1,242.97 | 26.9% |
| 11 | 281.70 | 1,184.74 | - | 5,457.62 | 968.24 | 16.2% |
| 12 | 159.08 | 838.80 | 29.2% | 6,280.94 | 761.72 | 12.2% |
| 13 | 161.89 | 665.17 | 14.7% | 3,816.97 | 635.02 | 7.5% |
| 14 | 129.03 | 543.06 | 10.3% | 2,533.62 | 531.91 | 6.1% |
| 15 | 118.63 | 466.04 | 6.5% | 1,132.61 | 452.16 | 4.7% |
| 16 | 128.28 | 407.95 | 4.9% | 1,099.09 | 394.06 | 3.4% |
| 17 | 102.99 | 355.73 | 4.4% | 892.94 | 351.78 | 2.5% |
| 18 | 121.27 | 318.89 | 3.1% | 1,094.77 | 317.42 | 2.0% |
| 19 | 114.33 | 288.22 | 2.6% | 1,397.27 | 287.68 | 1.7% |
| 20 | 120.39 | 263.92 | 2.1% | 1,062.88 | 262.42 | 1.5% |
| 21 | 95.47 | 240.89 | 1.9% | 790.06 | 238.87 | 1.4% |
| 22 | 98.91 | 221.44 | 1.6% | 641.12 | 219.51 | 1.1% |
| 23 | 100.45 | 206.66 | 1.2% | 486.45 | 203.57 | 0.9% |
| 24 | 84.78 | 193.87 | 1.1% | 348.49 | 191.03 | 0.7% |
| 25 | 72.41 | 185.36 | 0.7% | 766.11 | 182.60 | 0.5% |

Table 2. Computational time (sec.) and disutility function values for the uniform scenario.

| | Uniform Scenario | | | | | |
|---|---|---|---|---|---|---|
| | $T_{max} = 30$ | | | $T_{max} = 50$ | | |
| $p$ | CPU time (sec.) | Disutility | Marginal Benefit | CPU time (sec.) | Disutility | Marginal Benefit |
| 1 | 49.88 | Infeasible | - | 194.72 | Infeasible | - |
| 2 | 50.10 | Infeasible | - | 892.65 | Infeasible | - |
| 3 | 56.44 | Infeasible | - | 889.15 | Infeasible | - |
| 4 | 63.35 | Infeasible | - | 659.46 | Infeasible | - |
| 5 | 64.62 | Infeasible | - | 8,151.19 | Infeasible | - |
| 6 | 84.84 | Infeasible | - | 9,127.18 | Infeasible | - |
| 7 | 112.57 | Infeasible | - | 10,813.81 | Infeasible | - |
| 8 | 141.72 | Infeasible | - | 10,817.98 | Infeasible | - |
| 9 | 147.35 | Infeasible | - | 10,809.44 | 697.16 | - |
| 10 | 207.27 | Infeasible | - | 10,814.30 | 519.17 | 25.5% |
| 11 | 475.85 | 491.87 | - | 5,086.66 | 388.47 | 18.7% |
| 12 | 380.00 | 338.69 | 31.1% | 2,891.96 | 301.43 | 12.5% |
| 13 | 292.27 | 259.00 | 16.2% | 2,225.33 | 240.07 | 8.8% |
| 14 | 213.59 | 206.14 | 10.7% | 1,623.99 | 195.64 | 6.4% |
| 15 | 147.27 | 167.80 | 7.8% | 1,369.91 | 160.69 | 5.0% |
| 16 | 196.02 | 144.00 | 4.8% | 1,372.20 | 136.47 | 3.5% |
| 17 | 154.05 | 122.33 | 4.4% | 1,232.79 | 117.61 | 2.7% |
| 18 | 126.07 | 104.97 | 3.5% | 929.58 | 102.44 | 2.2% |
| 19 | 124.67 | 94.17 | 2.2% | 847.32 | 92.33 | 1.5% |
| 20 | 118.97 | 84.32 | 2.0% | 921.58 | 83.64 | 1.2% |
| 21 | 116.37 | 76.09 | 1.7% | 689.23 | 75.28 | 1.2% |
| 22 | 114.82 | 69.94 | 1.2% | 711.08 | 68.95 | 0.9% |
| 23 | 101.52 | 64.11 | 1.2% | 573.13 | 63.07 | 0.8% |
| 24 | 91.36 | 59.18 | 1.0% | 406.21 | 58.45 | 0.7% |
| 25 | 72.92 | 55.50 | 0.7% | 1,165.64 | 54.61 | 0.6% |

*Observation 3: If the number of drones is sufficiently large, employing more expensive, long-range drones does not add much benefit.*

We can compare the disutility values attained with the two drone ranges for various $p$, in order to obtain information for strategic decisions on which type of drones should be acquired. Table 3 shows the improvements between the disutility value obtained when $T_{max} = 50$ compared with $T_{max} = 30$. It is interesting that such an improvement is high for $p = 11$ (the minimum value for which $T_{max} = 30$ is feasible), namely 18.27% and 21.02% for the non-uniform and uniform scenarios, respectively; however, it rapidly falls and becomes almost negligible for larger fleet sizes (i.e., $p \geq 17$). This can be explained considering that, as the number of drones grows, each platform serves demands that are closer.

*Observation 4: If the fleet operator gradually increases the number p of drones, optimal locations are usually robust, that is, often the new location set includes the earlier location set.*

Table 4 shows the number of demand points assigned to each platform as $p$ increases from 11 to 25 (if a site is not selected, the entry is 0) for the non-uniform scenario when $T_{max} = 30$, and Table 5 shows the results for the $T_{max} = 50$ case. Observation 4 holds when $T_{max} = 30$ (Table 4), with two exceptions; platform at site #15 is optimal when $p = 11$, and then for $p = 12$ platform at site # 19 replaces it, and the other exception is when $p = 11$ and 12, platform at site # 7 is optimal, which is replaced by platform at site #16 when $p > 12$. Similar considerations hold for $T_{max} = 50$ (Table 5), for which the only exception is for platform at site #6 which is optimal for $p = 18$ and 19, it becomes non-optimal when $p = 21$, and it is optimal again for $p \geq 22$. Similar tables and observations hold also for the uniform scenario, for which the corresponding results are included in the Appendix (Table A1 and Table A2).

Table 3. Comparison between disutility values with $T_{max} = 30$ and $T_{max} = 50$ for both the scenarios.

| p | Non-Uniform Scenario | | Uniform Scenario | |
|---|---|---|---|---|
| | Difference | Improvement (%) | Difference | Improvement (%) |
| 11 | 216.50 | 18.27% | 103.39 | 21.02% |
| 12 | 77.08 | 9.19% | 37.26 | 11.00% |
| 13 | 30.15 | 4.53% | 18.92 | 7.31% |
| 14 | 11.15 | 2.05% | 10.49 | 5.09% |
| 15 | 13.88 | 2.98% | 7.11 | 4.24% |
| 16 | 13.89 | 3.41% | 7.53 | 5.23% |
| 17 | 3.95 | 1.11% | 4.71 | 3.85% |
| 18 | 1.47 | 0.46% | 2.53 | 2.41% |
| 19 | 0.54 | 0.19% | 1.84 | 1.96% |
| 20 | 1.50 | 0.57% | 0.68 | 0.80% |
| 21 | 2.02 | 0.84% | 0.81 | 1.07% |
| 22 | 1.93 | 0.87% | 0.99 | 1.41% |
| 23 | 3.10 | 1.50% | 1.04 | 1.63% |
| 24 | 2.84 | 1.47% | 0.72 | 1.22% |
| 25 | 2.76 | 1.49% | 0.89 | 1.61% |

Table 4. Number of demand points served by each platform for the non-uniform scenario when $T_{max} = 30$.

| Platform ID | P | | | | | | | | | | | | | | |
|---|---|---|---|---|---|---|---|---|---|---|---|---|---|---|---|
| | 11 | 12 | 13 | 14 | 15 | 16 | 17 | 18 | 19 | 20 | 21 | 22 | 23 | 24 | 25 |
| #1 | 10 | 8 | 8 | 7 | 6 | 6 | 6 | 6 | 5 | 5 | 5 | 5 | 5 | 5 | 4 |
| #2 | 10 | 9 | 8 | 7 | 7 | 6 | 6 | 5 | 6 | 6 | 6 | 6 | 5 | 5 | 5 |
| #3 | 0 | 0 | 0 | 0 | 0 | 0 | 0 | 0 | 0 | 0 | 0 | 0 | 0 | 4 | 4 |
| #4 | 9 | 9 | 8 | 7 | 8 | 7 | 6 | 6 | 6 | 5 | 5 | 4 | 4 | 5 | 4 |
| #5 | 9 | 9 | 8 | 8 | 7 | 7 | 6 | 7 | 7 | 6 | 6 | 6 | 6 | 5 | 4 |
| #6 | 7 | 7 | 6 | 6 | 5 | 5 | 4 | 4 | 4 | 4 | 4 | 4 | 4 | 3 | 3 |
| #7 | 10 | 8 | 0 | 0 | 0 | 0 | 0 | 0 | 0 | 4 | 4 | 4 | 4 | 4 | 4 |

| Platform ID | | | | | | | | | | | | | | | |
|---|---|---|---|---|---|---|---|---|---|---|---|---|---|---|---|
| #8  | 0  | 8  | 7  | 6 | 6 | 6 | 6 | 5 | 5 | 5 | 4 | 4 | 4 | 4 | 4 |
| #9  | 7  | 7  | 7  | 6 | 6 | 5 | 5 | 5 | 5 | 5 | 4 | 4 | 4 | 4 | 4 |
| #10 | 0  | 0  | 0  | 0 | 0 | 0 | 0 | 0 | 0 | 0 | 4 | 4 | 4 | 4 | 4 |
| #11 | 11 | 9  | 8  | 7 | 7 | 6 | 6 | 6 | 6 | 6 | 5 | 5 | 5 | 5 | 5 |
| #12 | 0  | 0  | 0  | 0 | 6 | 5 | 5 | 5 | 5 | 4 | 4 | 4 | 4 | 4 | 4 |
| #13 | 0  | 0  | 0  | 0 | 0 | 6 | 6 | 6 | 5 | 5 | 5 | 5 | 5 | 5 | 4 |
| #14 | 0  | 0  | 0  | 0 | 0 | 0 | 0 | 0 | 0 | 0 | 0 | 0 | 3 | 3 | 3 |
| #15 | 10 | 0  | 0  | 0 | 0 | 0 | 0 | 6 | 5 | 5 | 5 | 5 | 5 | 5 | 5 |
| #16 | 0  | 0  | 9  | 9 | 8 | 8 | 7 | 6 | 6 | 6 | 6 | 5 | 5 | 5 | 5 |
| #17 | 0  | 0  | 8  | 7 | 7 | 7 | 6 | 6 | 5 | 5 | 5 | 4 | 4 | 4 | 4 |
| #18 | 0  | 0  | 0  | 0 | 0 | 0 | 0 | 0 | 0 | 0 | 0 | 0 | 0 | 0 | 3 |
| #19 | 0  | 10 | 8  | 8 | 7 | 6 | 6 | 5 | 5 | 5 | 5 | 5 | 5 | 5 | 5 |
| #20 | 0  | 0  | 0  | 0 | 0 | 0 | 6 | 5 | 5 | 6 | 5 | 5 | 5 | 5 | 4 |
| #21 | 8  | 8  | 7  | 7 | 6 | 6 | 6 | 6 | 6 | 5 | 5 | 5 | 5 | 3 | 3 |
| #22 | 0  | 0  | 0  | 7 | 7 | 7 | 7 | 6 | 6 | 5 | 5 | 5 | 4 | 4 | 5 |
| #23 | 9  | 8  | 8  | 8 | 7 | 7 | 6 | 5 | 5 | 5 | 5 | 5 | 4 | 4 | 4 |
| #24 | 0  | 0  | 0  | 0 | 0 | 0 | 0 | 0 | 3 | 3 | 3 | 3 | 3 | 3 | 3 |
| #25 | 0  | 0  | 0  | 0 | 0 | 0 | 0 | 0 | 0 | 0 | 0 | 3 | 3 | 2 | 3 |

Table 5. Number of demand points served by each platform for the non-uniform scenario when $T_{max} = 50$.

| | p | | | | | | | | | | | | | | |
|---|---|---|---|---|---|---|---|---|---|---|---|---|---|---|---|
| Platform ID | 11 | 12 | 13 | 14 | 15 | 16 | 17 | 18 | 19 | 20 | 21 | 22 | 23 | 24 | 25 |
| #1  | 9  | 8 | 7 | 6 | 6 | 6 | 6 | 5 | 6 | 5 | 5 | 5 | 5 | 5 | 5 |
| #2  | 9  | 8 | 8 | 7 | 7 | 6 | 6 | 6 | 6 | 6 | 6 | 5 | 5 | 5 | 5 |
| #3  | 0  | 0 | 0 | 0 | 0 | 0 | 0 | 0 | 0 | 0 | 0 | 0 | 0 | 3 | 3 |
| #4  | 9  | 9 | 8 | 8 | 7 | 7 | 7 | 6 | 6 | 6 | 5 | 5 | 5 | 5 | 4 |
| #5  | 9  | 9 | 8 | 8 | 6 | 6 | 6 | 6 | 6 | 6 | 6 | 5 | 5 | 5 | 5 |
| #6  | 0  | 0 | 0 | 0 | 0 | 0 | 0 | 0 | 4 | 4 | 0 | 3 | 3 | 3 | 3 |
| #7  | 0  | 0 | 0 | 0 | 0 | 0 | 0 | 0 | 0 | 0 | 4 | 4 | 3 | 3 | 3 |
| #8  | 8  | 7 | 7 | 6 | 6 | 6 | 5 | 5 | 4 | 4 | 4 | 4 | 3 | 4 | 4 |
| #9  | 0  | 0 | 0 | 6 | 6 | 5 | 5 | 5 | 5 | 5 | 4 | 4 | 4 | 4 | 4 |
| #10 | 0  | 0 | 0 | 0 | 0 | 0 | 0 | 0 | 4 | 4 | 4 | 4 | 4 | 4 | 4 |
| #11 | 10 | 8 | 8 | 7 | 7 | 6 | 6 | 6 | 6 | 5 | 5 | 5 | 5 | 5 | 5 |
| #12 | 0  | 0 | 6 | 6 | 6 | 5 | 5 | 5 | 5 | 5 | 4 | 4 | 4 | 4 | 4 |
| #13 | 0  | 0 | 0 | 0 | 0 | 6 | 6 | 5 | 5 | 5 | 5 | 5 | 5 | 5 | 4 |
| #14 | 0  | 0 | 0 | 0 | 0 | 0 | 0 | 0 | 0 | 0 | 0 | 0 | 3 | 3 | 2 |
| #15 | 0  | 0 | 0 | 0 | 0 | 0 | 6 | 5 | 5 | 5 | 5 | 5 | 5 | 5 | 5 |
| #16 | 10 | 9 | 9 | 8 | 8 | 8 | 7 | 6 | 6 | 6 | 6 | 6 | 5 | 5 | 5 |
| #17 | 9  | 8 | 8 | 8 | 7 | 7 | 6 | 6 | 5 | 5 | 6 | 5 | 5 | 5 | 4 |
| #18 | 0  | 0 | 0 | 0 | 0 | 0 | 0 | 0 | 0 | 0 | 0 | 0 | 0 | 0 | 3 |
| #19 | 9  | 9 | 8 | 8 | 7 | 6 | 6 | 6 | 5 | 5 | 5 | 5 | 5 | 5 | 5 |
| #20 | 0  | 0 | 0 | 0 | 6 | 6 | 6 | 6 | 5 | 5 | 5 | 5 | 5 | 4 | 4 |
| #21 | 8 | 8 | 7 | 7 | 6 | 6 | 6 | 6 | 6 | 6 | 5 | 5 | 5 | 4 | 4 |
| #22 | 0 | 8 | 8 | 7 | 7 | 7 | 6 | 6 | 6 | 5 | 5 | 5 | 5 | 4 | 5 |
| #23 | 10 | 9 | 8 | 8 | 8 | 7 | 5 | 6 | 6 | 5 | 5 | 5 | 5 | 4 | 4 |
| #24 | 0 | 0 | 0 | 0 | 0 | 0 | 0 | 4 | 3 | 3 | 3 | 3 | 3 | 3 | 3 |
| #25 | 0 | 0 | 0 | 0 | 0 | 0 | 0 | 0 | 0 | 0 | 3 | 3 | 3 | 3 | 3 |

*3.2.2.2 Effects of using different time slot lengths.*

We performed a set of experiments to explore how different time slot lengths would affect the overall performance of model (1)-(7). We compared the results of the model obtained with 15-minute time slots ($w = 15$) versus the results obtained with 5-minute time slots ($w = 5$). In both the models, we considered a set up time of 15 minutes, equivalent to one time slot for the 15-minute time slot model, and to three time slots for the 5-minute time slot model. Figure 5 (and Table A3 in the Appendix) shows the results for short-range drones (i.e., $T_{max} = 30$) and non-uniform disutility functions; similar results and conclusions (not reported here) were obtained for long-range drones (i.e., $T_{max} = 50$). For different values of the parameter $p$, Figure 5 shows that: (i) the percentage of the time horizon needed to serve all the demand points, referred to as *makespan* in the primary y-axis (we used this metric to compare the quality of the solutions since changing the time slot duration causes the total disutility values and last used time slot to have slightly different values); and (ii) the computational time required to solve the models (secondary y-axis).

The main conclusions we derived from this set of experiments are:
- *The model with 5-minute time slots is computationally more difficult to solve.*

This is an expected result since the number of decision variables increases, indeed the number of time slots is 144 when $w = 5$, while it is 48 when $w = 15$. As it is shown in Figure 5, the computational time needed to solve the 5-minute time slots model reaches the time limit of three hours for most of the scenarios ($p = 11 - 22$, and $p = 25$ marked by asterisks), while it is usually less than two minutes when solving the model considering 15-minute time slots.
- *The model with 5-minute time slots shows a small improvement in the overall performance of the system compared with the model with 15-minute time slots.*

The percentage gap in the makespan obtained when solving the model with $w = 15$ vs $w = 5$ is computed as $\frac{Makespan\ (w=5)}{Number\ of\ timeslots\ (w=5)} - \frac{Makespan\ (w=15)}{Number\ of\ timeslots\ (w=15)}$ and it is always less than 5% (detailed values are reported in Table A3 in the Appendix). The fact that in some cases ($p = 12,13,15,16,17,18$) the reported makespan is higher for $w = 5$ is a consequence of the higher computational complexity. In fact, the comparison is made between the *best-found* solutions of the two models when the time limit was reached.

From these results we can conclude the following:

*Observation 5: Using time slots smaller than 15 minutes does not result in a significantly better makespan. On the other hand, it significantly increases the computational time required to solve the model.*

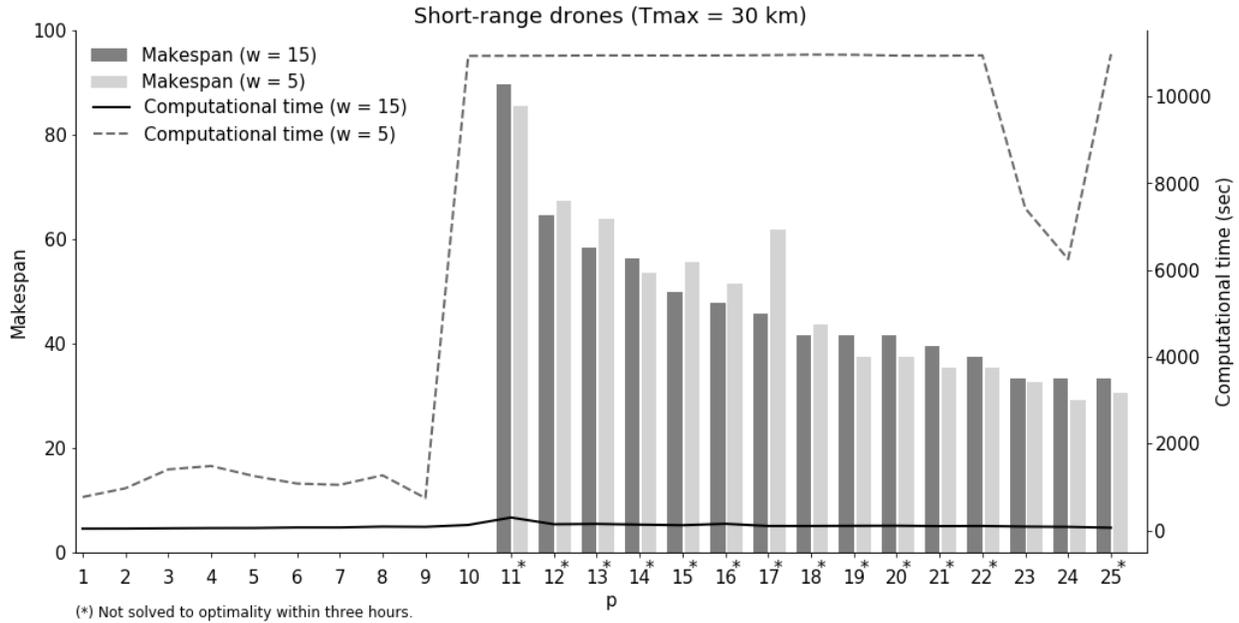

**Figure 5.** Comparison of makespan and computational times when solving the single-period model with 5-minute time slots ($w = 5$) and 15-minutes time slots ($w = 15$) for short-range drones.

*3.2.2.3 Effects of allowing up to two drones per platform.*

We now report the results of the experiments performed to analyze the effects of allowing more than one drone per platform. Specifically, we compare the results of model (1) - (7) when at most one drone is allowed per platform and when at most two drones are allowed per platform. To allow the use of up to two drones per platform, no modification of the structure of the mathematical model is needed, only the instances must be suitably specified. Specifically, we *duplicate* each candidate platform, so that there are now 50 possible potential location sites, and model (1) - (7) is solved on this new instance. If both the copies of the same site are selected, then that location is assigned two drones. For these experiments, parameter $p$ in the model corresponds to the fleet size, i.e., the number of used drones.

Figure 6 (and Table A4 in the Appendix) compares the results of the two versions of the model when either one drone is allowed per platform (Model-1) or up to two drones are allowed per platform (Model-2) for a fleet of short-range drone and non-uniform disutility functions. Similar results for long-range drones are reported in the Appendix in Figure A1 and Table A5. In order to compare the two versions of the model, we are showing the results for a number of available drones from 1 to 25 (even if the model would allow using up to 50 drones). Figure 6 shows the value of the objective function (primary y-axis) and the number of utilized platforms (secondary y-axis), versus the number of used drones (x-axis).

*Observation 6. With a small-medium fleet of drones, allowing more drones per platform is not beneficial.*

When the fleet size is less than or equal to 15 drones, the optimum solution is always to use one drone per platform. Hence, the solutions of the two models are the same.

*Observation 7. With a large fleet of drones, allowing more than one drone per platform reduces the number of required platforms.*

For a larger fleet (i.e., number of drones greater than 15) allowing multiple drones per platform reduces the number of utilized platforms and slightly improves the overall value of the total disutility. As shown in Figure 6, allowing up to two drones per platform can save up to 7 platforms (see the results corresponding to 24 used drones). Moreover, the larger the fleet, the larger is the improvement of the disutility function with a maximum improvement up to 14% (see Table A4 in the Appendix when the number of available drones is equal to 25).

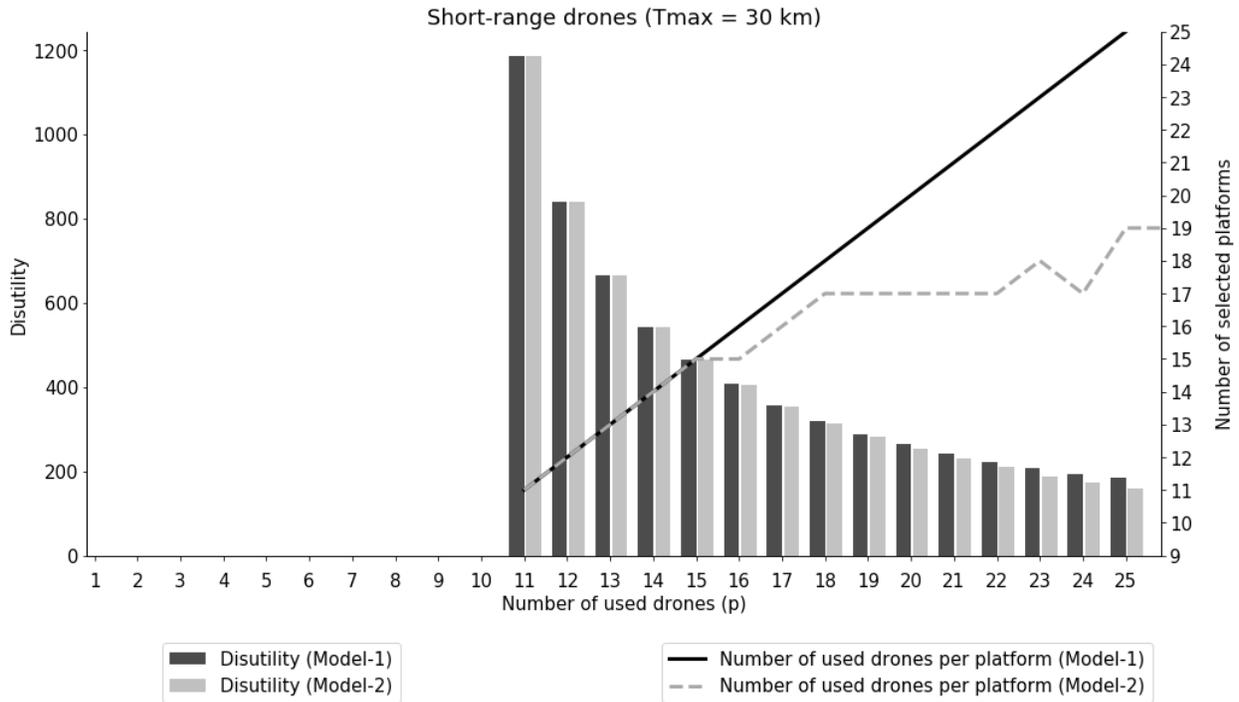

**Figure 6. Comparison of total disutility and of number of utilized platforms between Model-1 and Model-2 for short-range drones.**

*3.2.2.4 Effects of using a homogenous fleet of drones versus using a mixed fleet of drones.*

To address the effect of using a mixed fleet versus the use of a homogeneous fleet of drones, we compared the results obtained by running model (1) - (7) with those obtained by solving model (27)- (38) provided in the Appendix. In model (27)-(38), for a given value of $p$, it is assumed half of the drones are short-range ($N^{(1)} = \lceil \frac{p}{2} \rceil$) and half of the drones are long-range ones ($N^{(2)} = \lfloor \frac{p}{2} \rfloor$). We considered both the cases when one drone is allowed per platform, and the case where up to two drones are allowed per platform. Figure 7 (and Table A6 in the Appendix) provides a comparison of the disutility values obtained when adopting either only short-range drones, or only long-range drones, or a mix of drone types, and one drone per platform is allowed. Figure 8 (and Table A7 in the Appendix) provide similar results when up to two drones are allowed per platform. For both the cases we considered non-uniform disutility functions.

When only one drone per platform is allowed, the optimum disutility value for a mixed fleet of drones lies between the optimal values obtained when the fleet is composed of only short-range drones and when the fleet is composed of only long-range drones. Note that model (27)-(38) could not be solved to optimality within three hours on some of the scenarios (these are denoted with an asterisk in Figures 7 and 8). Let us take a closer look at the scenario with $p = 10$. When 10 long-range drones are available, all the demand points can be served with a total disutility of 1242.96 (Table A6 in the Appendix), while if a mixed fleet is available (i.e., five long-range and five short range drones) all the demand points can be served with a total disutility equal to 1405.53. Hence, although both cases are feasible, the total disutility value obtained using only long-range drones is lower than the total disutility value obtained with a mixed fleet of drones. The difference in total disutility decreases with larger fleet of drones, Table A6 in the Appendix provides the detailed results and shows the percentage gap improvement of using a fleet of long-range drones versus a mixed fleet of drones. The benefit is negligible (less than 1%) when the size of the fleet is greater than 14 drones, while the % gap is no more than 12% for a smaller fleet.

Hence, assuming drones with different range have different costs (even if cost of the drones is not explicitly accounted for in the model formulation), employing a mix of half short and half long-range drones when the fleet size is more than 13 drones maybe more cost-effective. The following observation can be derived:

*Observation 8: Assuming drones with different range have different costs, employing a mix of half short and half long-range drones for a large fleet could be more cost-effective than using a homogeneous fleet of drones.*

We observe a similar behavior (Figure 8 below and Table A7 in the Appendix) when up to two drones per platform are allowed. When adopting only short-range drones the smallest $p$ for which feasibility is achieved is $p = 11$; while, when adopting either only long-range drones or a mixed fleet, feasibility can be obtained for $p = 9$ and $p = 10$. Again, we observed that the disutility values obtained considering a mixed fleet of drones lie between the disutility values obtained with homogenous fleets of drones; it is worth mentioning that most of the models with a mixed fleet are not solved to optimality within three hours. Table A7 in the Appendix lists the disutility savings when only long-range drones are employed compared to the employment of a mixed fleet. Notwithstanding that some of the cases are not solved to optimality, the highest gain is 8.78%.

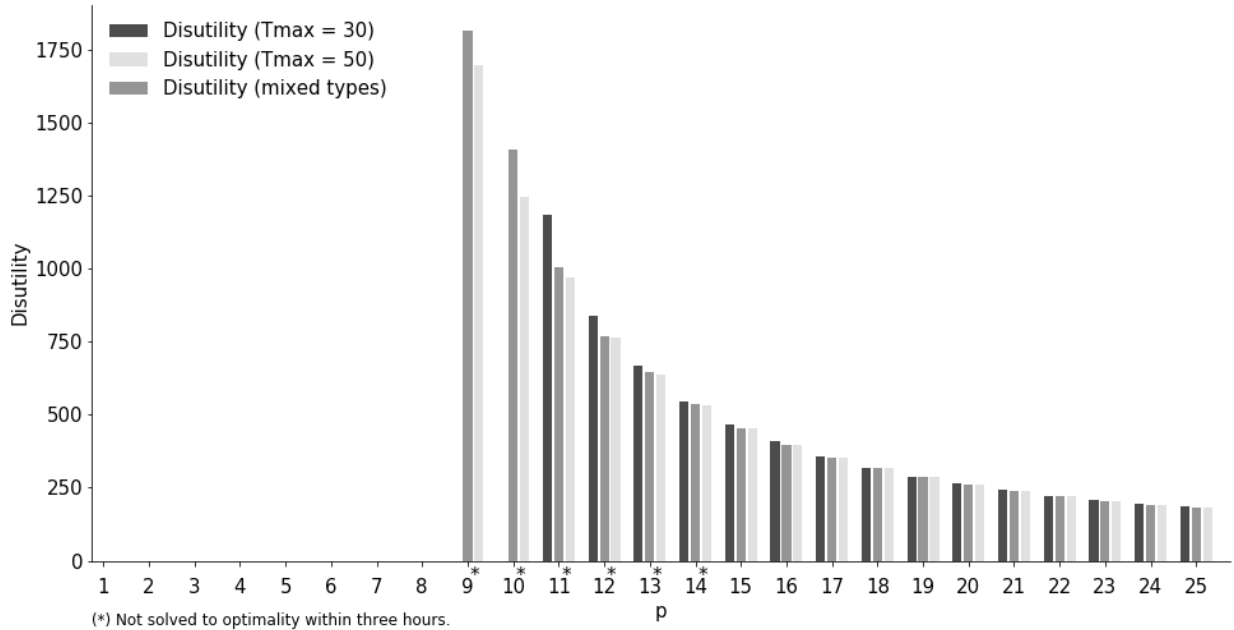

**Figure 7. Comparison of disutility values when using a homogenous fleet of drones vs a mixed fleet of drones, in the case that at most one drone is allowed per platform.**

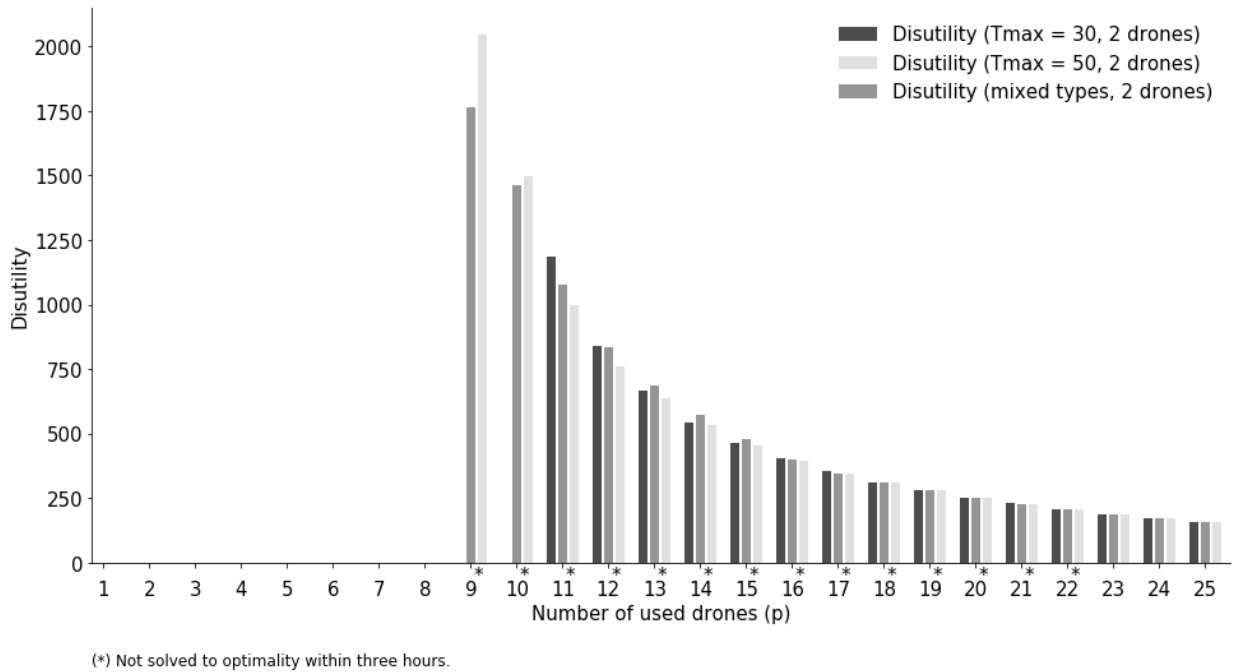

**Figure 8. Comparison of disutility values when using a homogenous fleet of drones vs a mixed fleet of drones, in the case that at most two drones are allowed per platform.**

# 4 The Two-period Problem

## 4.1 A time-slot formulation for the two-period DLS

In what follows we consider the two-period DLS, where we assume that the planning horizon is *two days*, and that during the night some platforms may be relocated, that is, each platform remains in the same site for one day, but some can be moved during the night. (That is, it is assumed that relocations, if any, can take place at night.) Also, it is assumed there may be a limit on the number of platforms that can be relocated.

A formulation of the two-period DLS is obtained by extending formulation (1) - (7) (a similar extension can be applied to formulation (27) - (38) in the Appendix). The meaning of variables $x_{itk}$ is the same as for the single-period formulation, except that now $t$ spans two days. Variables $y_i$ are replaced by $y_i^1$ and $y_i^2$ indicating when platform at site $i$ is used on day 1 and day 2, respectively. In order to count the number of platform relocations, we introduce a set of new variables $q_i, i \in P$, where $q_i = 1$ if site $i$ is used *exactly in only one of the two days*. Notice that although, in principle, one can think of moving a platform from site $i$ to site $i'$ and another platform from site $i'$ to site $i''$, such a solution is equivalent to a solution in which the platform at site $i'$ remains at $i'$ and the platform at site $i$ is moved to $i''$. Hence, w.l.o.g., we assume that a platform will only move overnight to a site which was empty during the first day, and the total number of platforms moved will be given by *half* the number of variables $q_i$ equal to 1. For example, if a platform is moved from $i$ to $j$, both $q_i = 1$ and $q_j = 1$ hold, but a single relocation takes place. In our formulation, $T_1$ and $T_2$ denote the sets of time slots of the two days respectively, and $T = T_1 \cup T_2$. We number the time slots $1, 2, \ldots, |T_1|, |T_1| + 1, |T_1| + 2, \ldots, |T_1| + |T_2|$.

$$\min \sum_{i \in P} \sum_{k \in D(i)} \sum_{t \in T} f_k\left(t - \frac{p_{ik}}{2}\right) x_{itk} \qquad (11)$$

$$\sum_{i \in P} \sum_{t \in T} x_{itk} = 1 \qquad \forall k \in D \qquad (12)$$

$$\left(\sum_{\substack{h \in D \\ h \neq k}} \sum_{\tau = t - p_{ik}}^{t} x_{i\tau k}\right) \leq M_{ik}(1 - x_{itk}) \qquad \forall i \in P, \forall t \in T_1, \forall k \in D(i) \qquad (13)$$

$$\left(\sum_{\substack{h \in D \\ h \neq k}} \sum_{\tau = t - p_{ik}}^{t} x_{i\tau k}\right) \leq M_{ik}(1 - x_{itk}) \qquad \forall i \in P, \forall t \in T_2, \forall k \in D(i) \qquad (14)$$

$$x_{itk} \leq y_i^1 \qquad \forall i \in P, \forall t \in T_1, \forall k \in D \qquad (15)$$

$$x_{itk} \leq y_i^2 \qquad \forall i \in P, \forall t \in T_2, \forall k \in D \qquad (16)$$

$$\sum_{i \in P} y_i^1 = p \qquad (17)$$

$$\sum_{i \in P} y_i^2 = p \qquad (18)$$

$$\sum_{t \in T_1 : t \leq p_{ik}} x_{itk} = 0 \qquad \forall i \in P, \forall k \in D \qquad (19)$$

$$\sum_{t \in T_2: t \leq 48 + p_{ik}} x_{itk} = 0 \qquad \forall i \in P, \forall k \in D \qquad (20)$$

$$q_i \geq y_i^1 - y_i^2 \qquad \forall i \in P \qquad (21)$$

$$q_i \geq y_i^2 - y_i^1 \qquad \forall i \in P \qquad (22)$$

$$\sum_{i \in P} q_i \leq 2Q \qquad (23)$$

$$x_{itk} \in \{0,1\} \qquad \forall i \in P, \forall t \in T, \forall k \in D \qquad (24)$$

$$y_i \in \{0,1\} \qquad \forall i \in P \qquad (25)$$

$$q_i \in \{0,1\} \qquad \forall i \in P \qquad (26)$$

Note that constraint sets (13) and (14), as well as (15) and (16), (17) and (18), and (19) and (20), refer to the two days separately. Constraints (21)-(23) account for the relocations, where parameter $Q \leq p$ specifies the maximum number of platform relocations that are allowed in this scenario.

### 4.2 Computational experiments for the two-period DLS

The computational results for the two-period model are discussed below. In all the experiments, the demand points and the platform sites are the same as for the single-period model and 15-minute time slots are considered. Also, the perishability and importance coefficients are the same as in the non-uniform scenario. The disutility function now spans two days. It is defined similarly as for the single-period model (i.e., equations (8)-(10)), however, for time slots 49, 50, ..., 96 a fixed penalty value in the second term of the disutility function given by equation (10) is added to represent overnight delay effects. Such a penalty accounts for the overnight break and it is equal to 50 if the due date of the demand point falls in the second day (i.e., $d_k \geq 49$), while it is equal to 100 otherwise.

We ran three different sets of experiments.

- *Short-due-dates experiments.* In these instances, the due dates of the demand points are exactly the same as in the single-period experiments (i.e., $d_k \leq 48, \forall k \in D$), but the time horizon now spans two days (96 time slots), that is, the drones may fly also on the second day to meet unfulfilled demands. Hence, some instances which were infeasible for the single-period model are now feasible.
- *Extended-due-dates experiments.* In these instances, the due dates of the demand points span both days (i.e., $d_k \leq 96, \forall k \in D$), and they are generated uniformly from time slot 1 to time slot 96.
- *Fewer-platforms experiments.* The disutility functions are the same as in the extended due dates experiments, but one site has been excluded from consideration. The aim of these experiments is to evaluate how sensitive total disutility is to the availability of potential platform sites.

All the experiments were run on a server with CPU Intel Xeon E5-2650 v3 @ 2.30GHz, with 128 GB of RAM, using CPLEX 12.8. A maximum time limit of three hours was set for each run.

The input data and the detailed results are available upon request to the authors. The details of the analysis are discussed in the remainder of the section. Before examining the computational results, we illustrate a sample solution. Figures 9a and 9b show the solution of the two-period model obtained when solving an instance of the extended-due-date experiments for the case $T_{max} = 30$, $p = 6$ and $Q = 4$ (i.e., four relocations are allowed). Figures 9a shows the schedule for day 1, while Figure 9b shows the schedule for day 2. In both the figures, black squares indicate unselected platform sites, red stars indicate selected platform sites used only on Day 1, blue pentagons indicate selected platform sites used only on Day 2, and purple hexagons indicate selected platform sites used on both days. Green circles are demand points served on day one, while grey circles are demand points served on day 2. Each selected platform site is linked to the allocated demand points. The number in any demand point is the time slot when the drone is back to the platform after serving the demand point. The solution shows that 74 demand points are served on day 1, while 26 demand points are served on day 2. The selected platforms sites in day 1 are (refer to Figure 4) #7, #11, #15, #20, #21 and #22. The selected platform sites in day 2 are #4, #6, #9, #11, #15, and #24. Hence, platform sites #11 and #15 are used in both days, while platform sites #7, #20, #21 and #22 are used in day 1, and then relocated at sites #4, #6, #9, and #24. Note that the model does not account for the cost of relocations, hence it does not specify which platform is relocated where. If necessary, a cost-effective relocation pattern may be defined by finding a minimum cost matching between platform sites active on day 1 and platform sites active on day 2.

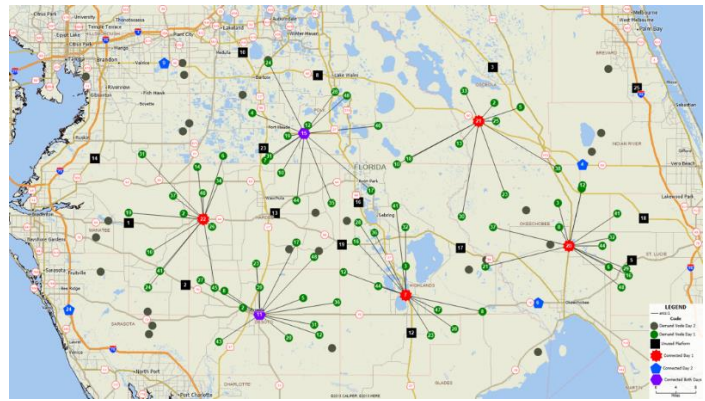

(a) Schedule on day 1.

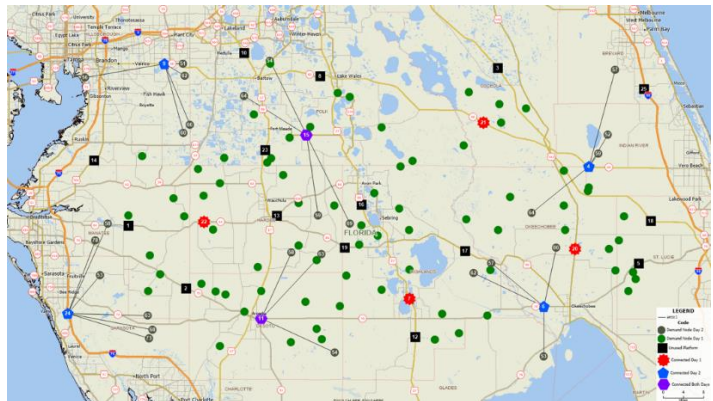

(b) Schedule on day 2.

Figure 9. example solution of the two-period model for one of the instances in the extended-due-date experiments when $T_{max} = 30$, $p = 6$, and $Q = 4$.

### 4.2.1 Short-due-dates experiments

Table 6 and Table 7 report the results of the two-period model for the short due-date experiments, that is when the due dates associated with each demand point are all in day 1. Several insights can be gained from examining the results, as discussed below.

- *Observation 9:* If some deliveries can be carried out on the second day, then fewer drones are necessary to obtain a feasible schedule.
- *Observation 10:* Most often, the marginal benefit of relocating one platform overnight is smaller than the marginal benefit of adding a drone.
- *Observation 11:* The marginal benefit of adding one drone is larger when the fleet is small and long-range drones are utilized.
- *Observation 12:* Under some specific circumstances, relocating a platform can be almost as effective as having one additional drone.

Table 6 shows the results for relevant $(p, Q)$ pairs when $T_{max} = 30$. For each $(p, Q)$ pair the table shows the value of the disutility function, the computational time (in seconds), the total number of demand points served on day 1, and the total number of demand points served on day 2. For each value of $p$, we report the results for those values of $Q$ for which a feasible solution was found. Specifically, the solutions shown in this table are also optimal. For values of $p \geq 11$, the value of the solution when $Q \geq 0$ coincides with the solution of the corresponding single-period model and results are therefore not shown here. Rows in grey denote the instances in which one additional relocation would not result in any improvement for that value of $p$, although CPU times sometimes differ. Table 7 shows the results when $T_{max} = 50$. For values of $p \geq 9$ the value of the solution when $Q \geq 0$ coincides with the solution of the corresponding single-period model and results are therefore not shown here. For each value of $p$, we show the results for those values of $Q$ for which a feasible solution was found. In no cases was CPLEX able to find the optimal solution within three hours CPU time limit imposed in these experiments, so for relevant $(p, Q)$ pairs, this table shows the value of the best two-period solution found and the gap with respect to the current lower bound returned by the solver.

Recall from Section 3.2.2.1 that when $T_{max} = 30$, the smallest value of $p$ to achieve feasibility in a single day is $p = 11$, while (Table 6) when deliveries can be performed on the second day a feasible schedule can be obtained with $p = 6$ by allowing $Q = 4$ relocations. Similarly, when $T_{max} = 50$, the smallest value of $p$ to achieve feasibility in a single day is $p = 9$. However (see Table 7), when deliveries can be performed on the second day, a feasible schedule can be obtained with $p = 5$ and $Q = 4$ relocations. Hence, considering an extended delivery period of two days, fewer drones are necessary to achieve a feasible schedule, if relocations are allowed (*Observation 9*). It is worth mentioning that the results reported in Table 6 and Table 7 for the cases when no relocation is allowed (that is, $Q = 0$) correspond also to the results we would obtain when running the single-period model with an extended time horizon of 24 hours (corresponding to 96 time slots). In these cases, the smallest value of $p$ to achieve feasibility in a single day would be $p = 10$ for $T_{max} = 30$, and $p = 6$ for $T_{max} = 50$.

Tables 6 and 7 also allow evaluating the tradeoff between the number of drones and the number of relocations. For example, consider Table 6, and compare the value of the disutility attained over two days when $p = 10$ with that obtained with the single-period model when $p = 11$. If we allow deliveries on the second day but no relocations are possible (i.e., $Q = 0$), then $p = 10$ drones are enough for fulfilling deliveries, and the total disutility is equal to 2,357.09. If we decide to add one drone and not to perform relocations the total disutility is equal to 1,184.74 (last row of Table 6, the same figure also appears in Table 1). This corresponds to a relative disutility improvement of 49.7% (i.e., $\frac{2,357.09-1,184.743}{2,357.09}$). If, instead of adding a drone, one decides to relocate one platform overnight (i.e., $p = 10$ and $Q = 1$), the corresponding attained disutility is 1,688.11, yielding an improvement equal to 28.4% (i.e., $\frac{2,357.09-1,688.11}{2,357.09}$). Hence, adding one drone is more beneficial that relocating one platform overnight (*Observation 10*).

Table 6. Results of the Short Due Dates Experiments for $T_{max} = 30$.

| $p$ | $Q$ | Disutility | CPU time (sec.) | # served on day 1 | # served on day 2 |
|---|---|---|---|---|---|
| 6 | 4 | 10,423.26 | 3,771.69 | 74 | 26 |
| 6 | 5 | 9,218.90 | 6,321.75 | 79 | 21 |
| 6 | 6 | 8,180.48 | 1,116.47 | 78 | 22 |
| 7 | 3 | 7,164.55 | 5,085.71 | 80 | 20 |
| 7 | 4 | 6,179.92 | 1,157.52 | 86 | 14 |
| 7 | 5 | 5,745.12 | 2,346.61 | 86 | 14 |
| 7 | 6 | 5,683.06 | 767.43 | 86 | 14 |
| 7 | 7 | 5,668.69 | 750.32 | 86 | 14 |
| 8 | 2 | 5,124.77 | 6,300.06 | 88 | 12 |
| 8 | 3 | 4,219.83 | 3,144.52 | 94 | 6 |
| 8 | 4 | 3,880.42 | 669.95 | 94 | 6 |
| 8 | 5 | 3,877.60 | 630.98 | 94 | 6 |
| 8 | 6 | 3,877.60 | 544.56 | 94 | 6 |
| 8 | 7 | 3,877.60 | 633.61 | 94 | 6 |
| 8 | 8 | 3,822.23 | 617.47 | 94 | 6 |
| 9 | 1 | 3,315.99 | 10,849.21 | 95 | 5 |
| 9 | 2 | 2,723.56 | 3,743.11 | 96 | 4 |
| 9 | 3 | 2,571.91 | 1,247.06 | 96 | 4 |
| 9 | 4 | 2,569.09 | 865.71 | 96 | 4 |
| 10 | 0 | 2,357.59 | 356.07 | 97 | 3 |
| 10 | 1 | 1,688.11 | 592.10 | 99 | 1 |
| 11 | 0 | 1,184.74 | 360.79 | 100 | 0 |

Table 7. Results of the Short Due Dates Experiments for $T_{max} = 50$.

| $p$ | $Q$ | Gap | Disutility (best solution found) | #served on day 1 | #served on day 2 |
|---|---|---|---|---|---|
| 5 | 4 | 4.22% | 12,563.08 | 64 | 36 |
| 5 | 5 | 1.67% | 12,187.40 | 64 | 36 |
| 6 | 0 | 18.09% | 9,470.00 | 74 | 26 |
| 6 | 1 | 3.38% | 7,915.45 | 77 | 23 |
| 6 | 2 | 2.77% | 7,783.82 | 77 | 23 |
| 6 | 3 | 15.92% | 7,783.82 | 77 | 23 |
| 6 | 4 | 3.01% | 7,685.32 | 77 | 23 |
| 6 | 5 | 3.35% | 7,685.32 | 77 | 23 |
| 6 | 6 | 3.64% | 7,685.32 | 77 | 23 |
| 7 | 0 | 23.48% | 5,780.49 | 88 | 12 |
| 7 | 1 | 8.48% | 4,794.93 | 90 | 10 |
| 7 | 2 | 3.94% | 4,643.94 | 90 | 10 |
| 7 | 3 | 5.95% | 4,643.94 | 90 | 10 |
| 7 | 4 | 6.40% | 4,643.94 | 90 | 10 |
| 7 | 5 | 3.64% | 4,621.82 | 90 | 10 |
| 7 | 6 | 3.69% | 4,621.82 | 90 | 10 |
| 7 | 7 | 2.74% | 4,615.61 | 90 | 10 |
| 8 | 0 | 6.41% | 2,511.17 | 99 | 1 |
| 8 | 1 | 5.31% | 2,480.88 | 99 | 1 |
| 8 | 2 | 5.59% | 2,480.88 | 99 | 1 |
| 8 | 3 | 4.81% | 2,480.88 | 99 | 1 |
| 8 | 4 | 3.63% | 2,476.18 | 99 | 1 |
| 8 | 5 | 5.20% | 2,476.18 | 99 | 1 |
| 8 | 6 | 5.38% | 2,476.18 | 99 | 1 |
| 8 | 7 | 3.85% | 2,476.18 | 99 | 1 |
| 8 | 8 | 5.39% | 2,476.18 | 99 | 1 |
| 9 | 0 | 0.01% | 1,688.72 | 100 | 0 |

Similar conclusions can be derived when $T_{max} = 50$. In this case the difference between the benefit of adding a drone versus the relocation of one platform is even larger. Specifically (see Table 7), when $p = 6$ and $Q = 0$, the disutility is equal to 9,470. If one relocation is allowed (i.e., $Q = 1$), the disutility equals 7,915.45, corresponding to an improvement of 16.4%. If a drone is added (i.e., $p = 7$ and $Q = 0$), the disutility equals 5,780.49, corresponding to an even greater improvement of 38.9%. Such a difference between the two benefit increases is consistent for each value of $p$ as shown in Table 8, where the relative benefit of relocating one platform is compared to that of adding one drone for $p = 6, 7$ and 8 (note that for $p \geq 9$, relocation does not help since all the demands can be served in a single day). Hence, these results show that adding one drone instead of relocating one platform is preferable with a small fleet (i.e., $p \leq 8$) of long-range drones (*Observation 11*). The explanation is that the number of deliveries carried out on the second day is not sensitive to the number of relocations allowed (unlike the case for $T_{max} = 30$), so the benefit of covering a significant number of deliveries on the first day, from having an additional drone, is much higher than optimizing the deliveries of the second day (which has a relatively larger disutility value anyhow).

**Table 8. Relative benefit of relocating one platform and relative benefit of adding one drone when $T_{max} = 50$.**

| p | Benefit of relocation | Benefit of one additional drone |
|---|---|---|
| 6 | 16.42% | 38.96% |
| 7 | 17.05% | 56.56% |
| 8 | 1.21% | 32.75% |

Finally, it is worth examining the scenario $p = 10$ when $T_{max} = 30$ (Table 6). A bit surprisingly, when $p = 10$ and $Q = 0$, three deliveries are scheduled in the second day in the optimal solution, and only one in the optimal solution when $p = 10$ and $Q = 1$. This may appear odd. In fact, for each item delivered on the second day, the fixed disutility penalty is added. So, since there is the possibility of delivering 99 items out of 100 on the first day, it is somewhat surprising that this does not occur when $Q = 0$. A careful analysis of these two solutions provides the following reason. There is one demand point (#87) which can only be served by platform at site #6, all other platform locations being farther than 30km. So, when $Q = 0$, one platform is constrained to be located at site #6, which is a very non beneficial choice and forces three deliveries to be scheduled on the second day. On the other hand, if a single relocation is allowed, one platform can be moved to site #6 overnight, and, hence, used to serve the single demand point #87. This example shows that in some circumstances, relocating a platform can be almost as effective as having one additional drone (*Observation 12*).

### 4.2.2 Extended-due-dates experiments

In this set of instances, the due dates of the demand points are uniformly distributed over two days. The main findings from this set of experiments can be summarized as follows:

- *Observation 13:* Relocation is beneficial when the fleet size is small. When a large fleet is available, then relocation may not beneficial.
- *Observation 14*: The benefit of relocation is more evident with a small fleet of short-range drones, since this allows finding feasible schedules which would not be possible otherwise.
- *Observation 15*: When relocation is effective, a small number of relocations suffices to achieve the total benefit in terms of disutility. This relocation number decreases when the size of the fleet increases.
- *Observation 16*: A small fleet of long-range drones is more effective when relocation is allowed.

Table 9 shows the results for relevant $(p, Q)$ pairs when $T_{max} = 30$ and $T_{max} = 50$, respectively. For each $(p, Q)$ pair the table shows the value of the disutility function and the computational time (in seconds). We report the results for those values of $p$ for which a feasible solution was found for at least one value of $Q$. If the time limit was not reached, the displayed solution is also optimum. As before, rows in grey denote the instances for which one additional relocation would not produce any improvement for that value of $p$, although CPU times might differ.

These results show well the tradeoff between the number of drones and the number of relocations required to get feasible solutions. When $T_{max} = 30$ (see Table 9), if no relocations are allowed, at least $p = 10$ drones are needed to satisfy all the demands. If fewer drones are available, relocations can partially compensate. Indeed, a feasible solution can be obtained with $p = 6$ if at least $Q = 4$ relocations are allowed. When $T_{max} = 50$ (see Table 9), a feasible solution is obtained when $p = 5$ and using three relocations (i.e., $Q = 3$); relocations are not needed for $p \geq 9$. Hence, relocation is useful when the fleet is small, i.e., $p \leq 10$ (*Observation 13*), and more so for short-range drones (i.e., $T_{max} = 30$), in which case relocation allows achieving a feasible schedule, which would not be possible otherwise (*Observation 14*).

Observe from Table 9 (i.e., when $T_{max} = 30$), the values of $p$ for which relocation is beneficial. There is a minimum number $Q_{min}$ of relocations that give a feasible solution (i.e., $Q_{min} = 4, 3, 2, 1, 0$, for $p = 6, 7, 8, 9, 10$, respectively). Relocations do not provide any benefit for $p \geq 12$. From Table 9 (i.e., when $T_{max} = 50$); when $p = 5$ we have $Q_{min} = 3$, and when $p = 6, 7, 8$ we have $Q_{min} = 0$. So, overall, when relocations are effective, a small number of them are needed to achieve the total benefit in terms of disutility. This number decreases when the size of the fleet increases (*Observation 15*).

The comparisons between the results obtained with different values of $T_{max}$ are summarized in Table 10 which shows the difference in the objective function values for $T_{max} = 30$ and for $T_{max} = 50$ and the corresponding relative improvements. These results should be compared with those shown in Table 3 for the single-day model. As for the single-day model, observe that the largest differences occur for small values of $p$ and $Q$ (improvements greater than or equal to 20% are outlined in grey in Table 10). Hence, the employment of more expensive long-range drones is useful when the fleet of drones is small (*Observation 3*). Additionally, if relocation is allowed, a small fleet of long-range drones becomes more effective than a small fleet of short-range drones (*Observation 16*). This is a reasonable finding, since range limitations can be compensated by having more drones and/or more relocations. Specifically, when $p = 6$ (and $Q = 4$ relocations are allowed), using drones with $T_{max} = 50$ brings a 29% saving. This saving is 33% for $p = 7$ (and $Q = 3$), 45% for $p = 8$ (and $Q = 2$), 43% for $p = 9$ (and $Q = 1$), 43% for $p = 10$ (with no relocations), 19% for $p = 11$ (with no relocations) and then it steadily decreases until it becomes less than 1% for $p \geq 17$.

Finally, from the computational-time viewpoint, observe from Table 9 that when $T_{max} = 30$, the hardest instance was the one with $p = 9$ and $Q = 1$, the only one for which the 3-hr. CPU time limit was reached. However, the optimality gap was less than 1%. For the case $T_{max} = 50$, the model appears significantly harder to solve. In fact, the time limit was reached in most cases with small values of $p$ and $Q$, with a gap occasionally as large as 19%.

Table 9 . Results of the Extended Due Dates experiments

| $T_{max} = 30$ | | | | $T_{max} = 50$ | | | |
|---|---|---|---|---|---|---|---|
| $p$ | $Q$ | Disutility | CPU time (sec.) | $p$ | $Q$ | Disutility | CPU time (sec.) |
| 6 | 4 | 7,706.25 | 2,941.15 | 5 | 3 | 10,705.36 | > 3 hrs |
| 6 | 5 | 6,658.10 | 3,340.87 | 5 | 4 | 9,083.74 | > 3 hrs |
| 6 | 6 | 5,738.45 | 1,688.28 | 5 | 5 | 9,077.58 | > 3 hrs |
| 7 | 3 | 5,148.78 | 2,631.79 | 6 | 0 | 6,331.03 | > 3 hrs |
| 7 | 4 | 4,445.40 | 1,772.89 | 6 | 1 | 6,017.90 | > 3 hrs |
| 7 | 5 | 4,142.24 | 1,471.41 | 6 | 2 | 5,573.84 | > 3 hrs |
| 7 | 6 | 4,093.43 | 576.34 | 6 | 3 | 5,529.95 | > 3 hrs |
| 7 | 7 | 4,080.27 | 2,013.00 | 6 | 4 | 5,499.59 | > 3 hrs |
| 8 | 2 | 3,763.11 | 3,519.94 | 6 | 5 | 5,499.59 | > 3 hrs |
| 8 | 3 | 3,088.98 | 3,193.03 | 6 | 6 | 5,490.50 | > 3 hrs |
| 8 | 4 | 2,961.63 | 793.07 | 7 | 0 | 3,493.35 | > 3 hrs |
| 8 | 5 | 2,955.49 | 940.74 | 7 | 1 | 3,454.99 | > 3 hrs |

| 8 | 6 | 2,915.32 | 583.85 | 7 | 2 | 3,446.02 | > 3 hrs |
| 9 | 1 | 2,621.17 | > 3 hrs | 7 | 3 | 3,442.86 | > 3 hrs |
| 9 | 2 | 2,122.97 | 1,127.11 | 7 | 4 | 3,436.78 | > 3 hrs |
| 9 | 3 | 2,048.14 | 1,247.40 | 7 | 5 | 3,434.48 | > 3 hrs |
| 9 | 4 | 2,046.72 | 1,037.00 | 7 | 6 | 3,430.80 | > 3 hrs |
| 10 | 0 | 2,002.32 | 476.99 | 8 | 0 | 2,099.80 | > 3 hrs |
| 10 | 1 | 1,474.40 | 700.20 | 8 | 1 | 2,052.87 | > 3 hrs |
| 11 | 0 | 1,114.31 | 690.44 | 9 | 0 | 1,485.14 | 9,861.04 |
| 11 | 1 | 1,062.18 | 266.89 | 10 | 0 | 1,147.36 | 7,289.12 |
| 12 | 0 | 809.16 | 279.03 | 11 | 0 | 907.99 | 5,849.54 |
| 13 | 0 | 649.79 | 531.16 | 12 | 0 | 727.61 | 2,971.35 |
| 14 | 0 | 535.01 | 211.79 | 13 | 0 | 612.06 | 2,515.09 |
| 15 | 0 | 461.51 | 252.20 | 14 | 0 | 523.64 | 2,296.14 |
| 16 | 0 | 405.61 | 248.29 | 15 | 0 | 447.77 | 1,709.68 |
| 17 | 0 | 354.51 | 234.36 | 16 | 0 | 391.30 | 1,389.25 |
| 18 | 0 | 318.44 | 236.24 | 17 | 0 | 350.07 | 1,077.45 |
| 19 | 0 | 288.03 | 447.68 | 18 | 0 | 315.99 | 1,440.42 |
| 20 | 0 | 263.79 | 215.20 | 19 | 0 | 286.55 | 1,583.45 |
| 21 | 0 | 240.75 | 202.91 | 20 | 0 | 261.65 | 1,207.17 |
| 22 | 0 | 221.31 | 207.70 | 21 | 0 | 238.36 | 1,244.59 |
| 23 | 0 | 206.58 | 193.43 | 22 | 0 | 219.43 | 971.89 |
| 24 | 0 | 193.79 | 162.01 | 23 | 0 | 203.48 | 805.53 |
| 25 | 0 | 185.29 | 152.68 | 24 | 0 | 190.94 | 712.01 |
| - | - | - | - | 25 | 0 | 182.57 | 1,042.30 |

Table 10. Comparison in the disutility function between the scenario with $T_{max} = 30$ and the scenario with $T_{max} = 50$.

| p | Q | Disutility $T_{max} = 30$ | Disutility $T_{max} = 50$ | Difference | Relative Improvement |
| --- | --- | --- | --- | --- | --- |
| 6 | 4 | 7,706.25 | 5,499.59 | 2,206.66 | 29% |
| 6 | 5 | 6,658.10 | 5,499.59 | 1,158.51 | 17% |
| 6 | 6 | 5,738.45 | 5,490.50 | 247.94 | 4% |
| 7 | 3 | 5,148.78 | 3,442.86 | 1,705.92 | 33% |
| 7 | 4 | 4,445.40 | 3,436.78 | 1,008.62 | 23% |
| 7 | 5 | 4,142.24 | 3,434.48 | 707.76 | 17% |
| 7 | 6 | 4,093.43 | 3,430.80 | 662.63 | 16% |
| 8 | 2 | 3,763.11 | 2,052.87 | 1,710.24 | 45% |
| 8 | 3 | 3,088.98 | 2,052.87 | 1,036.11 | 34% |
| 8 | 4 | 2,961.63 | 2,052.87 | 908.76 | 31% |
| 8 | 5 | 2,955.49 | 2,052.87 | 902.61 | 31% |
| 9 | 1 | 2,621.17 | 1,485.14 | 1,136.02 | 43% |
| 9 | 2 | 2,122.97 | 1,485.14 | 637.83 | 30% |
| 9 | 3 | 2,048.14 | 1,485.14 | 562.99 | 27% |
| 10 | 0 | 2,002.32 | 1,147.36 | 854.95 | 43% |
| 10 | 1 | 1,474.40 | 1,147.36 | 327.03 | 22% |
| 11 | 0 | 1,114.31 | 907.99 | 206.32 | 19% |
| 11 | 1 | 1,062.18 | 727.61 | 334.56 | 31% |
| 12 | 0 | 809.16 | 727.61 | 81.55 | 10% |
| 13 | 0 | 649.79 | 612.06 | 37.74 | 6% |
| 14 | 0 | 535.01 | 523.64 | 11.37 | 2% |
| 15 | 0 | 461.51 | 447.77 | 13.75 | 3% |
| 16 | 0 | 405.61 | 391.30 | 14.31 | 4% |
| 17 | 0 | 354.51 | 350.07 | 4.44 | 1% |
| 18 | 0 | 318.44 | 315.99 | 2.45 | 1% |
| 19 | 0 | 288.03 | 286.55 | 1.47 | 1% |
| 20 | 0 | 263.79 | 261.65 | 2.14 | 1% |
| 21 | 0 | 240.75 | 238.36 | 2.40 | 1% |
| 22 | 0 | 221.31 | 219.43 | 1.88 | 1% |
| 23 | 0 | 206.58 | 203.48 | 3.10 | 1% |
| 24 | 0 | 193.79 | 190.94 | 2.84 | 1% |
| 25 | 0 | 185.29 | 182.57 | 2.72 | 1% |

*4.2.2 Experiments with Fewer Platform Sites*

In these experiments we examined the impact of availability of fewer sites on the overall disutility. Results are shown in Table 11.

In particular, we re-ran the instances used for the extended-due-date experiments after removing one site from consideration among those sites that were selected, for all values of $p$; excluding more than one site led to infeasibility for any value of $p$ and $Q$ for the $T_{max} = 30$ case. For a given value of $p$, Table 11 report the percentage disutility increases when one site, specifically site #1 (see Figure 9), was excluded from consideration when $Q = 0$ for the cases $T_{max} = 30$ and $T_{max} = 50$, respectively. Observe that when $T_{max} = 30$ these percentages monotonically increases, reaching, for $p = 24$, 8.48% for $T_{max} = 30$ and 8.57% for $T_{max} = 50$. The reason for this behavior is that most likely, as long as $p$ is small, the excluded site can be easily replaced by some other location, that is, there are several platform configurations which are close to being equivalent. When $p$ is large, there are obviously fewer possible configurations to choose from (actually, the configuration becomes fixed when $p = 24$ and one site is excluded), so the percentage difference increases between the two cases, with and without excluding a site. Nonetheless, in absolute terms the difference between disutility values remains very small for $p \geq 18$.

For $T_{max} = 50$ and $Q = 0$ (see Table 11), the exclusion of one platform has a similar impact for all $p \geq 6$. When $p = 6,7,8$ the differences are larger because, we believe, the solver reached the CPU time limit on these instances and the solutions are not optimal.

Table 11. Results of the experiments with extended due dates and fewer platform sites (24 vs. 25) when $Q = 0$.

| | $T_{max} = 30$ | | $T_{max} = 50$ | |
|---|---|---|---|---|
| $p$ | Disutility difference | Disutility percentage increase | Disutility difference | Disutility percentage increase |
| 6 | - | - | 693.02 | 10.95% |
| 7 | - | - | 2,077.83 | 59.48% |
| 8 | - | - | 247.66 | 11.79% |
| 9 | - | - | 98.93 | 6.66% |
| 10 | 0.00 | 0.00% | 67.89 | 5.92% |
| 11 | 39.71 | 3.56% | 50.96 | 5.61% |
| 12 | 36.78 | 4.55% | 39.13 | 5.38% |
| 13 | 34.92 | 5.37% | 34.29 | 5.60% |
| 14 | 36.50 | 6.82% | 35.17 | 6.72% |
| 15 | 32.05 | 6.94% | 33.33 | 7.44% |
| 16 | 27.24 | 6.72% | 26.99 | 6.90% |
| 17 | 25.42 | 7.17% | 21.23 | 6.06% |
| 18 | 18.01 | 5.66% | 16.06 | 5.08% |
| 19 | 13.93 | 4.84% | 14.48 | 5.05% |
| 20 | 13.80 | 5.23% | 13.50 | 5.16% |
| 21 | 13.16 | 5.47% | 12.45 | 5.22% |
| 22 | 12.47 | 5.64% | 11.66 | 5.32% |
| 23 | 12.57 | 6.09% | 13.22 | 6.50% |
| 24 | 16.43 | 8.48% | 16.36 | 8.57% |

## 5 Summary and Future Research

### 5.1 Summary

Motivated by observations on the need of emergency delivery of medical products to demand points that are not well connected by a roadway network, we have introduced a drone delivery problem that hitherto has not been addressed. Each demand point may include a preferred due time (or due date) or may require the consideration of perishability of the delivered medical product, or both. A major realistic consideration in the problem is the limited flying range to serve each demand point from the location of the take-off platform from which the drone serves the demand point. Another major characteristic of our new problem is that the drones operate out of mobile platforms which may be, for example, vans that can move on usable roads from a current platform site to another potential platform site. The third major feature of the model is that, besides recommending location and relocation of platforms at sites, it schedules delivery based on costs or disutilities of individual demands.

In the problem addressed, we assumed that we are given a set $P$ of $m$ candidate platform sites and a set $D$ of $n$ demand points. A total of $p \leq m$ platforms must be located at the sites, where each platform hosts one drone. If a platform is at site $i$, a drone departing from it can serve any demand point located within a specified limited range from it. A drone can only carry one package at a time, and it must head back to the platform after each delivery.

The paper addressed two variants of this problem. In the first variant, we are given a fixed time, for example, one day, from morning to evening, to deliver the products demanded. In the second variant, we have two periods, say two days, to deliver the demanded products. In both these problem variants, the objective is to minimize the total disutility for the delivery of the required packages of products. In the second variant, it is assumed that the mobile platforms may be moved in the evening of the first day for second day operations. We first observed that when time is modeled with a continuous variable $t$, the resulting mathematical formulation cannot be solved for realistic instances in a reasonable amount of time. To come up with a method that provides useful results, we developed an optimization model which essentially boils down to scheduling activities in *time-and-space slots*, where a "utility" is earned for timely delivery at a demand point in any given slot. In the model, we divided the operating period into a set $T$ of time slots. Associated with individual demand k,there is a non-decreasing disutility $f_k(t)$ function of time, that may be given by a continuous function of $t$ if it is specified, or by table of disutility values at the time slots corresponding to each value time $t$. We underscore that the strength of this approach can include both the preferred due time of the individual demand at point $k$ and the perishability of the product demanded.

Extensive computational experiments on the time-slot formulation demonstrated that, indeed, the results are available in reasonable computational times. Results for the single-period model included the following insights:

[1] a minimum number of platforms are needed for problem feasibility;
[2] this minimum number of platforms needed for problem feasibility decreases when the range of the drones increases;
[3] for feasible problems, the solution specifies (a) on which sites the given number $p$ platforms should be located, (b) from which platform which demand points should be served, and (c) the sequence in which the demand should be served from any given located platform;
[4] the marginal benefit of having one additional drone/platform monotonically strictly decreases;
[5] many of the sites are robust in the sense that as $p$ increases, they are still optimal; and
[6] time slot length should be carefully chosen since it significantly affects the computational efficiency of the model. A 15-minute time slot seems to be a reasonable choice.

For the two-period model, the problem solutions demonstrated, among several insights, the following:

[1] compared with the single-period case, fewer drones are necessary to achieve a feasible schedule when some deliveries can be carried out on the second day;
[2] relocation in period 2 is necessary when the fleet size is small;
[3] the minimum number of drones needed for problem feasibility decreases when the range of the drones is increased; and
[4] as for the single-period model, the two-period solution specifies (a) where the platforms should be located, and, (b) if necessary, where they should be relocated on the second day, and (c) which demands should be served in which sequence.

Finally, we emphasize that the developed formulations are flexible enough to easily accommodate the case when multiple drones can operate from each platform and can also be easily modified to allow for a mixed fleet of drones.

In summary, the major contributions of the paper are follows: (1) it addressed a new important drone-based delivery problem that solved for three sets of decisions through one optimization model, (i) the locations of drone platforms, (ii) the allocation of demands to each located platform, and (iii) the sequencing of the deliveries from the selected platform sites; (2) it developed a practical useful and flexible model for its solution; and (3) it conducted extensive computational results that demonstrated the model's usefulness.

### 5.2 Potential future applications and research

During emergency situations, where people are stranded without means of transport, there may be critical needs which need immediate attention, especially when the communication and power infrastructure may be damaged. Examples of such needs are (a) essential medicines, e.g., insulin and blood pressure pills, (b) minimum food/water for survival, (c) communication mechanism, e.g., satellite phones, (d) heat source (in case of very cold weather), and (e) fuel if it can be used for necessary functions, e.g., diesel fuel for electricity generation. Another example of a potential new application is by the U.S. Navy which often performs a wide range of remote activities off their vessels including battlefield support and international humanitarian operations. These activities include rescues at sea, transport of emergency personnel and relief supplies, emergency relief operations, delivery of medical supplies to troops, and many other activities. If drones are used for some of these activities, models like the ones developed in this paper could be useful.

This paper limited the drones to <u>delivery</u> of products. Drones could be used also for pick up. Logistic optimization models may be needed that provide an itinerary for each drone when the fleet is given a list of delivery or pick-up *chores* at specific locations; an example itinerary might include its take off point, its first drop off chore, then a pick-up chore, then the return to the base, then a pickup of another item from the inventory, and so on. Critical issues addressed by such models and methods would be (1) *Timeliness* of specific deliveries or pickups since there could be issues of perishability (e.g., medicine, insulin, blood samples etc.), deadlines and time windows for delivery/pickup, (2) *Limitations of each drone*, in particular its <u>range</u> (per charge if it is electrical or per fuel tank if it is powered by liquid fuel) and <u>payload</u>, and (3) *Substitutability* among delivered items, for examples, a high-end satellite phone may be substituted for low-end one; the last issue is especially pertinent when the inventory available to the fleet is limited. It would be extremely useful to develop logistical models for recommending actions for a fleet of drones that is provided a list of delivery or pick-up chores at specific locations.

The other enhancements of the models could include the followings:

(a) Only drones perform the delivery tasks. If the moveable platform itself could perform some tasks, then the problems become more complex and TSP/VRP considerations must be included in the models.
(b) It may be useful for the platforms to move while the associated drone is busy delivering a product. In this case some dynamics of the mobile platforms and their synchronization with the associated drones must be included in the optimization models.

## Acknowledgements

The authors acknowledge colleague Matthew A. Treglia for helping with the graphical visualizations of the solutions of the various problem instances, ASU graduate student Mr. Kiarash Ghasemzadeh for assisting in the literature review, and ASU graduate student Mr. Brandon Mathis for assisting in developing the Florida case study.